# Weak backward error analysis for Langevin process

**Marie Kopec**

**Abstract** We consider numerical approximations of stochastic Langevin equations by implicit methods. We show a weak backward error analysis result in the sense that the generator associated with the numerical solution coincides with the solution of a modified Kolmogorov equation up to high order terms with respect to the stepsize. This implies that every measure of the numerical scheme is close to a modified invariant measure obtained by asymptotic expansion. Moreover, we prove that, up to negligible terms, the dynamic associated with the implicit scheme considered is exponentially mixing.

## 1 Introduction

In the last decades, backward error analysis has become a powerful tool to analyze the long time behavior of numerical schemes applied to evolution equations (see [9,15,22]). The main idea can be described as follows: Let us consider an ordinary differential equation of the form

$$\dot{y}(t) = f(y(t)),$$

where $f : \mathbb{R}^d \to \mathbb{R}^d$ is a smooth vector field, and denote by $\phi_t^f(p)$ the associated flow. By definition, a numerical method defines for a small time step $\delta$ an approximation $\Phi_\delta$ of the exact flow $\phi_\delta^f$: We have for bounded $p \in \mathbb{R}^d$, $\Phi_\delta(p) = \phi_\delta^f(p) + \mathcal{O}(\delta^{r+1})$ where $r$ is the order of the method.

The idea of backward error analysis is to show that $\Phi_\delta$ can be interpreted as the exact flow $\phi_\delta^{f_\delta}$ of a modified vector field defined as a series in powers of $\delta$

$$f_\delta = f + \delta^r f_r + \delta^{r+1} f_{r+1} + ...,$$

where $f_l$, $l \geq r$ are vector fields depending on the numerical method. In general, the series defining $f_\delta$ does not converge, but it can be shown that for bounded $y$, we have for arbitrary $N$

$$\Phi_\delta(y) = \phi_\delta^{f_\delta^N}(y) + C_N \delta^{N+1},$$

M. Kopec
ENS Cachan Bretagne, campus de Ker Lann, avenue Robert Schumann, 35170 Bruz, France
E-mail: marie.kopec@bretagne.ens-cachan.fr

where $f_\delta^N$ is the truncated series:
$$f_\delta^N = f + \delta^r f_r + ... + \delta^N f_N.$$

Under some analytic assumptions, the constant $C_N \delta^{N+1}$ can be optimized in $N$, so that the error term in the previous equation can be made exponentially small with respect to $\delta$.

Such a result is very important and has many applications in the case where $f$ has some strong geometric properties, such as Hamiltonian (see [8,9,15,21,22]). In this situation, and under some compatibility conditions on the numerical method $\Phi_\delta$, the modified vector field $f_\delta$ inherits the structure of $f$.

More recently, these ideas have been extended in some situations to Hamiltonian PDEs : First in the linear case [5] and then in the semi linear case (nonlinear Schrödingier or wave equations), see [6,7].

In [4], the authors give a weak backward error analysis for SDEs defined on the d-dimensional torus. The aim of this article is to extend the result of [4] to the Langevin process on $\mathbb{R}^{2d}$. This process is defined by the stochastic Hamiltonian equation:

$$\begin{aligned} dq(t) &= M^{-1} p(t) dt, \\ dp(t) &= -\partial_q V(q(t)) dt - \gamma p(t) dt + \sigma M^{1/2} dW(t), \end{aligned} \quad (1.1)$$

where $V$ is the potential energy function of a classical model for a molecular system, $M$ is a mass matrix, $\gamma$ is a free parameter, the friction coefficient and the term $\sigma dW(t)$ is a fluctuation term bringing energy into the system. This equation can be used to give an approximation of the following integral

$$\int_{\mathbb{R}^{2d}} \phi(q,p) Z^{-1} \exp(-\beta H(q,p)) dq dp,$$

where $\beta = \frac{1}{k_B T} = \frac{2\gamma}{\sigma^2}$, T is temperature, $k_B$ is Boltzmann's constant, $H(q,p) = \frac{1}{2} p^t M^{-1} p + V(q)$ and $Z = \int_{\mathbb{R}^{2d}} \exp(-\beta H(q,p)) dq dp$ (see [3,13,16] for more explanations).

In this work, we investigate the weak error which concerns the law of the solution. Let us recall that given a SDE in $\mathbb{R}^d$ of the form
$$dX = f(X) dt + g(X) dW, \quad (1.2)$$

discretized by an explicit Euler scheme $(X_p)$ with time step $\delta$, then, under assumptions on $f$, $g$ and $\phi : \mathbb{R}^d \to \mathbb{R}$ (see [19,20,25,26]), the explicit Euler scheme $(X_p)$ has weak error order 1:

$$|\mathbb{E}(\phi(X_p)) - \mathbb{E}(\phi(X(p\delta)))| \leq c(\phi, T) \delta, \quad p = 0, ..., \lfloor T/\delta \rfloor, \quad T > 0.$$

Error estimates on long times for elliptic and hypoelliptic SDEs have already been proved, especially in the case of explicit scheme ([24,25,27]) or on the torus ([19]). In [24,25], it is shown that for a sufficiently small time step, the explicit Euler scheme defines an ergodic process and that the invariant measure of the Euler scheme is close to the invariant measure of the SDE. In [27], under the assumption of the existence of a unique invariant measure associated to the SDE, Talay and Tubaro have shown that the weak error and the invariant measure associated to the Euler scheme can be expanded in powers of the time step $\delta$. The assumptions on $f$ and $g$ used in [24,25,27] are restrictive. Moreover, the results describe in these papers are only for explicit schemes.
In [19], a larger class of schemes is studied. It is shown that given an elliptic or hypoelliptic SDE, the ergodic averages provided by a class of implicit and explicit schemes are asymptotically close to the average of the invariant measure of the SDE. The authors also show an expansion in expectation of the invariant measure for any time-step. Unfortunately, they work on the d-dimensional torus.

The behavior in long time of approximations of the Langevin equation have already been study, but the authors need that the time step is small enough [18,26] or work on the torus for the position [14]. In [18], under some assumptions, it is shown, in particular in the case of the Langevin equation, that for sufficiently small time step, two kinds of implicit schemes are ergodic processes. They also show that the invariant measures associated with theses schemes converge to the invariant measure of the Langevin equation. In [26], Talay shows the exponential convergence of the solution associated to the Kolmogorov equation and, for a sufficiently small time step $\delta$, an expansion with respect to $\delta$ of the invariant measure to the implicit scheme which is close to the invariant measure of the SDE. In [14], the authors study the Langevin equation and work with explicit splitting methods. They provide error estimates on the invariant distribution for small step size, and compare the sampling bias obtained for various choices of explicit splitting method.

In this paper, we work on $\mathbb{R}^{2d}$ and $V$ and all its derivatives are not necessarily bounded but have polynomial growth. The aim of this paper is, under assumptions on $V$, to show a weak backward error analysis result : We show an expansion with respect to the time step $\delta$ of $\mathbb{E}\phi(q_k, p_k)$ where $(q_k, p_k)$ is an implicit scheme. Unlike [18, 26], we do not need that the time step is small enough and our assumptions are less restrictive than in [26].

The idea to extend the backward error analysis to SDE has already be studied in [1,4,12,23]. In [1], the authors use this approach to construct new methods of weak order two to approximate stochastic differential equations. In [4], the authors study a SDE defined on the d-dimensional torus and its approximation by the explicit Euler scheme. They show, without restriction on the time step, an expansion of $\mathbb{E}\phi(X_k)$ where $X_k = (q_k, p_k)$ is the explicit Euler scheme. In [12], the author shows the same kind of result but for the overdamped Langevin equation and for implicit schemes. In [23], Shardlow consider SDE with additive noise ($g$ does not depend of $X$). He has shown that it is possible to build of a modified SDE associated with the Euler scheme, but only at the first step, i.e. for $N = 2$. In this case, he is able to write down a modified SDE:
$$d\tilde{X} = \tilde{f}(\tilde{X})dt + \tilde{g}(\tilde{X})dW,$$
such that
$$\left|\mathbb{E}\big(\phi(X_k)\big) - \mathbb{E}\big(\phi(\tilde{X}(k\delta))\big)\right| \leq c_1(\phi, T)\delta^2, \ k = 0, ..., \lfloor T/\delta \rfloor, \ T > 0.$$

In this paper, we take the approach describes in [4, 12]. We show that the generator associated with the process solution of the SDE coincides with the solution of a modified Kolmogorov equation up to high order terms with respect to the stepsize. It is known that given $\phi : \mathbb{R}^d \to \mathbb{R}$ and denoting by $(q_{q^0}(t), p_{p^0}(t))$ the solution of the SDE (1.1) satisfying $(q(0), p(0)) = (q^0, p^0)$, the function $u$ defined for $t \geq 0$ and $(q^0, p^0) \in \mathbb{R}^{2d}$ by $u(t, q^0, p^0) = \mathbb{E}(\phi(q_{q^0}(t), p_{p^0}(t)))$ satisfies the Kolmogorov equation
$$\partial_t u = Lu,$$
where $L$ is the Kolmogorov operator associated with the SDE.

We show that with the numerical solution, we can associate a modified Kolmogorov operator of the form
$$\mathscr{L}(\delta, q, p, \partial_q, \partial_p) = L(q, p, \partial_q, \partial_p) + \delta L_1(q, p, \partial_q, \partial_p) + \delta^2 L_2(q, p, \partial_q, \partial_p) + ... \quad ,$$
where $L_l$, $l \geq 1$ are some modified operators of order $2l + 2$. The series does not converge but we consider truncated series:
$$L^{(N)}(\delta, q, p, \partial_q, \partial_p) = L(q, p, \partial_q, \partial_p) + \delta L_1(q, p, \partial_q, \partial_p) + \delta^2 L_2(q, p, \partial_q, \partial_p) + ... + \delta^N L_N(q, p, \partial_q, \partial_p) \quad .$$

Note that this operator is no longer of order 2 and we can not define easily a solution to the modified equation
$$\partial_t v^N(t, q, p) = L^{(N)}(\delta, q, p, \partial_q, \partial_p) v^N(t, q, p).$$
However, in our case, we can build an approximated solution $v^{(N)}$ such that
$$\| \mathbb{E}(\phi(q_k, p_k)) - v^{(N)}(k\delta, .) \|_{\mathscr{C}} \leq c_2(\phi, N)\delta^N, \quad k = 0, ..., [T/\delta], \quad T > 0,$$

where $\mathscr{C}$ is an appropriate space. As the constant $c_2$ does not depend of $T$, we have an approximation result valid on very long times. We also show that there exists a modified invariant measure for $L^{(N)}(\delta,q,p,\partial_q,\partial_p)$.

The main tools are to find schemes who have moments of all order bounded in time, the exponential convergence to equilibrium of $u$ and all its derivatives and the hypo-ellipticity of the Poisson equation, i.e. the equation $L(q,p,\partial_q,\partial_p)u = h$. The second tool is also used in [19]. A proof of the exponential convergence to equilibrium of $u$ and all its derivatives can be found in [26], but we need to now more precisely the dependance of the bound. Moreover, using a result describes in [18], we simplify the proof in [26] (see Appendix for more details).

In section 2, we introduce the SDE and the assumptions that we need. We also introduce the numerical schemes that we use. Then, we give some results on the Kolmogorov operator and on the solution of the Kolmogorov equation. In section 3, we give an asymptotic expansion of the weak error. In section 4, we study the modified operator $\mathscr{L}$ and its approximation. In section 5, we analyze the long time behavior of $v^{(N)}$.

## 2 Preliminaries

2.1 Presentation of the SDE

We warn the reader that constants may vary from line to line during the proofs, and that in order to use lighter notations we usually forget to mention dependence of the parameters. We use the generic notation $C$ for such constants.

In all the article, we write the dot product of two vectors $q = (q_1,...,q_d) \in \mathbb{R}^d$ and $p = (p_1,...,p_d) \in \mathbb{R}^d$ as
$$\langle q,p \rangle = \sum_{i=1}^d q_i p_i.$$

For a multi-index $\mathbf{k} = (k_1,...,k_d) \in \mathbb{N}^d$, we set $|\mathbf{k}| = k_1 + ... + k_d$ and for a function $\phi \in C^\infty(\mathbb{R}^d)$
$$\partial_q^\mathbf{k} \phi(q) = \frac{\partial^{|\mathbf{k}|} \phi(q)}{\partial^{k_1} q_1 ... \partial^{k_d} q_d}, \quad q = (q_1,...,q_d) \in \mathbb{R}^d.$$

For $p \in \mathbb{R}^d$, we set $\partial_p = (\frac{\partial}{\partial p_1},...,\frac{\partial}{\partial p_d})^\top$ and for $k \in \mathbb{N}^*$, $\partial_p^k$ is the differential of order $k$. We also use the following notation
$$\mathscr{C}_{pol}^\infty(\mathbb{R}^d) = \{f \in C^\infty(\mathbb{R}^d) \text{ such that } f \text{ and all its derivatives have polynomial growth}\}$$
$$= \{f \in C^\infty(\mathbb{R}^d) \text{ such that for all } \mathbf{k} = (k_1,...,k_d) \in \mathbb{N}^d, \exists C_\mathbf{k}, n_\mathbf{k}$$
$$\text{such that for all } q \in \mathbb{R}^d, \quad |\partial_q^\mathbf{k} f(q)| \leq C_\mathbf{k}(1+|q|^{2n_\mathbf{k}})\}.$$

Let $(\Omega, \mathscr{F}, \mathscr{F}_t, \mathbf{P}), t \geq 0$, be a filtered probability space and $W(t) = (W_1(t),...,W_d(t))$ be a $d$-dimensional $\{\mathscr{F}_t\}_{t \geq 0}$-adapted standard Wiener process. We want to give a similar result on $\mathbb{R}^{2d}$ of [4] and [12] for a process $(q(t), p(t))_{t \in \mathbb{R}_+}$ on $\mathbb{R}^{2d}$ which verifies the stochastic Hamiltonian differential system
$$q(t) = q(0) + M^{-1} \int_0^t p(s) ds,$$
$$p(t) = p(0) - \int_0^t \partial_q V(q(s)) ds - \gamma \int_0^t p(s) ds + \sigma M^{1/2} W(t), \quad t > 0,$$

where $\gamma > 0$ is a friction coefficient, $\sigma > 0$ and $M$ is a positive diagonal mass matrix. The function $V : \mathbb{R}^d \to \mathbb{R}$, is $C^\infty$ and verifies the following conditions:

**B-1:** The function $V$ is semi-convex : There exist a bounded function $V_1 \in C^\infty(\mathbb{R}^d)$ with bounded derivatives and a convex function $V_2 \in C^\infty(\mathbb{R}^d)$ such that $V = V_1 + V_2$.

**B-2:** There exist $\kappa > 0$ and $\beta \in ]0,1[$ such that for all $q \in \mathbb{R}^d$

$$\frac{1}{2}\langle \partial_q V(q), q \rangle \geq \beta V(q) + \gamma^2 \frac{\beta(2-\beta)}{8(1-\beta)}|q|^2 - \kappa. \tag{2.1}$$

**B-3:** $V \in \mathscr{C}^\infty_{pol}(\mathbb{R}^d)$.

**B-4:** There exists a constant $\kappa_1$ such that for all $q \in \mathbb{R}^d$, $V(q) \geq \kappa_1|q|^2$.

**Remark 2.1** *We can replace assumption **B**-4 by*

$$\exists \kappa_1, \kappa_2 : \quad V(q) \geq \kappa_1|q|^2 - \kappa_2 \quad \text{for all } q \in \mathbb{R}^d,$$

*however it is more convenient to assume that $V \geq 0$.*

*A consequence of the semi-convexity of $V$ (assumption **B**-1) is that there exists a positive constant $\theta$ such that for all $q \in \mathbb{R}^d$ and $h \in \mathbb{R}^d$,*

$$\partial_q^2 V(q)(h, h) \geq -\theta|h|^2. \tag{2.2}$$

*In the following, $\theta$ will be called the constant of semi-convexity of $V$.*

*A consequence of the assumption **B**-2 is the dissipativity inequality: There exists a strictly positive real number $\beta_1$ such that*

$$\langle q, \partial_q V(q) \rangle \geq \beta_1|q|^2 - \kappa \quad \text{for all } q \in \mathbb{R}^d. \tag{2.3}$$

*Moreover, a polynomial $V$ growing at infinity like $|q|^{2k}$, where $k \in \mathbb{N}^*$, will satisfy the assumptions.*

To slightly simplify the presentation that follows, we make the change of variables $q \to M^{-1/2}q$, $p \to M^{1/2}p$, with a corresponding adjustment of the potential. This is equivalent to assuming $M = I$ and the new potential verifies same assumptions as $V$.

Hence, for now on, we consider the following equation:

$$\begin{aligned} q(t) &= q(0) + \int_0^t p(s)ds, \\ p(t) &= p(0) - \int_0^t \partial_q V(q(s))ds - \gamma \int_0^t p(s)ds + \sigma W(t), \quad t > 0, \end{aligned} \tag{2.4}$$

where $V$ satisfies the conditions **B**.

**Remark 2.2** *Under these assumptions, we have that $(q(t), p(t))$ is well-defined for all $t > 0$ (see Chapter III, Theorem 4.1 in [11]).*

We also consider the Hamiltonian function $H$ defined for $(q, p) \in \mathbb{R}^{2d}$ by

$$H(q, p) := \frac{1}{2}|p|^2 + V(q).$$

We define the Kolmogorov operator $L$ associated with (2.4) for $\phi \in C^\infty(\mathbb{R}^{2d})$ and $(q,p) \in \mathbb{R}^{2d}$ by

$$L\phi(q,p) := \langle p, \partial_q \phi(q,p) \rangle - \langle \partial_q V(q) + \gamma p, \partial_p \phi(q,p) \rangle + \frac{\sigma^2}{2}\sum_{i=1}^d \frac{\partial^2}{\partial p_i \partial p_i}\phi(q,p).$$

For the study of (2.4), it is useful to define the Lyapunov function $\Gamma$ defined for all $(q,p) \in \mathbb{R}^{2d}$ by

$$\Gamma(q,p) = \frac{1}{2}|p|^2 + V(q) + \frac{\gamma}{2}\langle q, p \rangle + \frac{\gamma^2}{4}|q|^2 + 1. \tag{2.5}$$

Under assumptions **B**, the function $\Gamma$ verifies the following properties:

**Lemma 2.3** *We have for all $(q,p) \in \mathbb{R}^{2d}$*

$$\Gamma(q,p) \geq \frac{1}{8}|p|^2 + \frac{\gamma^2}{12}|q|^2 + 1. \tag{2.6}$$

*Moreover, for every $\ell \geq 1$, there exist strictly positive real numbers $a_\ell$ and $d_\ell$ such that we have for every $(q,p) \in \mathbb{R}^{2d}$*

$$L\Gamma^\ell(q,p) \leq -a_\ell \Gamma^\ell(q,p) + d_l. \tag{2.7}$$

**Remark 2.4** *We can find a similar proof of Lemma 2.3 in [18].*

*Proof* We first show (2.6). We have, for all $\varepsilon > 0$ and $(q,p) \in \mathbb{R}^{2d}$,

$$\gamma \langle p,q \rangle \geq -\frac{\varepsilon}{2}|p|^2 - \frac{\gamma^2}{2\varepsilon}|q|^2,$$

then, using the positivity of $V$, we get for all $(q,p) \in \mathbb{R}^{2d}$ and $\varepsilon > 0$

$$\Gamma(q,p) \geq (\frac{1}{2} - \frac{\varepsilon}{4})|p|^2 + \gamma^2(\frac{1}{4} - \frac{1}{4\varepsilon})|q|^2 + 1.$$

Taking $\varepsilon = \frac{3}{2}$, we show (2.6).

We now show the other property of $\Gamma$ (2.7). We first do the case $\ell = 1$. We have for any $(q,p) \in \mathbb{R}^{2d}$ and $\varepsilon > 0$

$$\Gamma(q,p) \leq \frac{1}{2}|p|^2 + V(q) + \frac{\gamma}{2}(\frac{\varepsilon}{2}|p|^2 + \frac{1}{2\varepsilon}|q|^2) + \frac{\gamma^2}{4}|q|^2 + 1.$$

Multiplying by $\beta$, taking $\varepsilon = \frac{2(1-\beta)}{\gamma\beta}$ and using assumption **B-2** on the derivative of $V$, we get for all $(q,p) \in \mathbb{R}^{2d}$

$$\beta\Gamma(q,p) \leq \frac{1}{2}|p|^2 + \beta V(q) + \gamma^2 \frac{\beta(2-\beta)}{8(1-\beta)}|q|^2 + \beta \leq \frac{1}{2}|p|^2 + \frac{1}{2}\langle \partial_q V(q), q \rangle + \kappa + \beta. \tag{2.8}$$

Using this inequality, we can now show (2.7) for $\ell = 1$. Indeed, we have for all $(q,p) \in \mathbb{R}^{2d}$

$$L\Gamma(q,p) = -\frac{\gamma}{2}|p|^2 - \frac{\gamma}{2}\langle q, \partial_q V(q) \rangle + \frac{d\sigma^2}{2}$$

$$\leq \frac{d\sigma^2}{2} + (\kappa + \beta)\gamma - \beta\gamma\Gamma(q,p).$$

We have shown (2.7) for $\ell = 1$ with $d_1 = \frac{d\sigma^2}{2} + \gamma(\kappa + \beta)$ and $a_1 = \beta\gamma$.

Let $\ell \geq 2$ be an integer. We now calculate $L\Gamma^\ell$. We have, for $(q,p) \in \mathbb{R}^{2d}$ and $i \in \{1,...,2d\}$, with the notation $x = (q,p)$,

$$\frac{\partial}{\partial x_i}\Gamma^\ell(q,p) = \ell\Gamma^{\ell-1}(q,p)\frac{\partial}{\partial x_i}\Gamma(q,p),$$

$$\frac{\partial^2}{\partial x_i \partial x_i}\Gamma^\ell(q,p) = \ell\Gamma^{\ell-1}(q,p)\frac{\partial^2}{\partial x_i \partial x_i}\Gamma(q,p) + \ell(\ell-1)\Gamma^{\ell-2}(q,p)(\frac{\partial}{\partial x_i}\Gamma(q,p))^2,$$

then we get for all $(q,p) \in \mathbb{R}^{2d}$

$$L\Gamma^\ell(q,p) = \ell\Gamma^{\ell-1}(q,p)L\Gamma(q,p) + \frac{\ell(\ell-1)}{2}\sigma^2\Gamma^{\ell-2}(q,p)\left|p + \frac{\gamma}{2}q\right|^2.$$

Using for $(q,p) \in \mathbb{R}^{2d}$

$$\left|p+\frac{\gamma}{2}q\right|^2 \leq 2\Gamma(q,p)$$

and computations for $\ell = 1$, we get for $(q,p) \in \mathbb{R}^{2d}$

$$L\Gamma^\ell(q,p) \leq \ell\Gamma^{\ell-1}(q,p)\left(d_1 - a_1\Gamma(q,p)\right) + \ell(\ell-1)\sigma^2\Gamma^{\ell-1}(q,p)$$
$$= -a_1\ell\Gamma^\ell(q,p) + \ell(d_1+\ell-1)\Gamma^{\ell-1}(q,p).$$

Using the fact that for any $x \in \mathbb{R}$ and $\varepsilon > 0$

$$|x|^{\ell-1} \leq \varepsilon|x|^\ell + C_\varepsilon,$$

we get

$$L\Gamma^\ell(q,p) \leq (-a_1\ell + \varepsilon\ell(d_1+\ell-1))\Gamma^\ell(q,p) + \ell(d_1+\ell-1)C_\varepsilon.$$

Choosing $\varepsilon = \frac{a_1}{\ell(d_1+\ell-1)}$, $a_\ell = a_1(\ell-1)$ and $d_\ell$ large enough, we obtain (2.7) for $\ell$. □

Using assumptions **B** and properties of $\Gamma$, we have the following result on the moment of the solution:

**Lemma 2.5** *Let a process $(q(.), p(.))$ which satisfies (2.4). We have for each $\ell \geq 1$ such that $\mathbb{E}|q(0)|^{2\ell} < \infty$ and $\mathbb{E}|p(0)|^{2\ell} < \infty$ that there exist strictly positive real numbers $k_\ell$, $\alpha_\ell$ and $C_\ell$ such that for all $t > 0$*

$$\mathbb{E}\left(|q(t)|^{2\ell} + |p(t)|^{2\ell}\right) \leq C_\ell\left(1 + (\mathbb{E}|p(0)|^{k_\ell} + \mathbb{E}|q(0)|^{k_\ell})\exp(-\alpha_\ell t)\right).$$

*Proof* Let $N \in \mathbb{N}^*$ be fixed. We consider

$$\tau_N = \inf\{t, \text{ such that } \max(|q(t)|, |p(t)|) \geq N\}.$$

Let $\ell \in \mathbb{N}^*$ be fixed. We will show the following: There exists a positive real number $C$ such that for all $t \geq 0$

$$\mathbb{E}\left(\Gamma^\ell(q(t),p(t))\right) \leq C\left(1 + \Gamma^\ell(q(0),p(0))\exp(-\alpha t)\right). \tag{2.9}$$

Let $a_\ell$ be the constant defined in (2.7). Let $\alpha < a_\ell$ be fixed. We apply the Itô formula to $\Gamma^\ell(q(t), p(t))\exp(\alpha t)$. We obtain for all $t \geq 0$

$$\Gamma^\ell\big(q(t \wedge \tau_N), p(t \wedge \tau_N)\big)\exp\big(\alpha(t \wedge \tau_N)\big)$$
$$= \Gamma^\ell(q(0),p(0)) + \alpha\int_0^{t\wedge\tau_N}\Gamma^\ell(q(s),p(s))\exp(\alpha s)ds$$
$$+ \int_0^{t\wedge\tau_N}\langle p(s), \partial_q\Gamma^\ell(q(s),p(s))\rangle\exp(\alpha s)ds$$
$$- \int_0^{t\wedge\tau_N}\langle \partial_q V(q(s)) + \gamma p(s), \partial_p\Gamma^\ell(q(s),p(s))\rangle\exp(\alpha s)ds$$
$$+ \sigma\int_0^{t\wedge\tau_N}\partial_p\exp(\alpha s)\langle\Gamma^\ell(q(s),p(s)), dW(s)\rangle$$
$$+ \frac{\sigma^2}{2}\int_0^{t\wedge\tau_N}\sum_{i=1}^d\frac{\partial^2}{\partial p_i\partial p_i}\Gamma^\ell(q(s),p(s))\exp(\alpha s)ds$$
$$= \Gamma^\ell(q(0),p(0)) + \alpha\int_0^{t\wedge\tau_N}\Gamma^\ell(q(s),p(s))\exp(\alpha s)ds$$
$$+ \int_0^{t\wedge\tau_N}L\Gamma^\ell(q(s),p(s))\exp(\alpha s)ds$$
$$+ \sigma\int_0^{t\wedge\tau_N}\exp(\alpha s)\langle\partial_p\Gamma^\ell(q(s),p(s)), dW(s)\rangle.$$

The stochastic integral is a square integrable martingale because $\Gamma(q(.), p(.))$ is bounded on $[0, t \wedge \tau_N]$. Thus its average vanishes. Using the positivity of $\Gamma$, the inequality (2.7) and $\alpha < a_l$, we obtain for all $t \geq 0$

$$\mathbb{E}\Big(\Gamma^l\big(q(t \wedge \tau_N), p(t \wedge \tau_N)\big) \exp\big(\alpha(t \wedge \tau_N)\big)\Big) \leq \Gamma^l(q(0), p(0)) + \frac{d_l}{\alpha}\mathbb{E}\big(\exp(\alpha(t \wedge \tau_N))\big)$$
$$+ (\alpha - a_l)\mathbb{E}\Big(\int_0^{t \wedge \tau_N} \Gamma^l(q(s), p(s))\exp(\alpha s)ds\Big)$$
$$\leq \Gamma^l(q(0), p(0)) + \frac{d_l}{\alpha}\mathbb{E}\big(\exp(\alpha(t \wedge \tau_N))\big).$$

Using Fatou's lemma in the left hand side and Monotone convergence Theorem on the right hand side of the last inequality, we have (2.9). To conclude, we use the property (2.6) and polynomial growth of $\Gamma$.

**Remark 2.6** *We can found an other proof of Lemma 2.5 in [26].*

2.2 Numerical schemes

For a small time step $\delta$ and $x \in \mathbb{R}^{2d}$, the classical explicit Euler method applied to (1.2), is defined for $i = 1, ..., 2d$, by $X_0 = x$ and the formula

$$X_{n+1}^i = X_n^i + \delta f^i(X_n) + \sum_{\ell=1}^m g_\ell^i(X_n)(W^\ell((n+1)\delta) - W^\ell(n\delta)), \quad n \geq 0. \tag{2.10}$$

The ordinary Euler scheme (2.10) may be unstable when the coefficients of the differential equation (2.4) are unbounded (see [18]). We are led to avoid explicit schemes. In fact, we study two different implicit schemes. For a small time step $\delta > 0$, we consider an implicit split-step scheme defined by $q_0 = q(0)$, $p_0 = p(0)$ and for $n \in \mathbb{N}^*$

$$\begin{cases} q_{n+1} = q_n + \delta p_n^* \\ p_n^* = p_n - \delta \gamma p_n^* - \delta \partial_q V(q_{n+1}), \\ p_{n+1} = p_n^* + \sigma(W((n+1)\delta) - W(n\delta)) = p_n^* + \sqrt{\delta}\sigma \eta_n, \end{cases} \tag{2.11}$$

where $\eta_n = (\eta_{n,1}, ..., \eta_{n,d})$ is a $\mathbb{R}^d$-valued random variable and $\{\eta_{n,i} : n \in \mathbb{N}, i \in \{1, d\}\}$ is a collection of i.i.d. real-valued random variables satisfying $\eta_{1,1} \sim \mathcal{N}(0, 1)$. We also consider the implicit Euler scheme defined by $q_0 = q(0)$, $p_0 = p(0)$ and for $n \in \mathbb{N}^*$

$$\begin{cases} q_{n+1} = q_n + \delta p_{n+1} \\ p_{n+1} = p_n - \delta \partial_q V(q_{n+1}) - \gamma \delta p_{n+1} + \sqrt{\delta}\sigma \eta_n, \end{cases} \tag{2.12}$$

where $\eta_n = (\eta_{n,1}, ..., \eta_{n,d})$ is as above.

Using the following Lemma, we get that these two schemes are well-defined for $\delta < \delta_0 := \frac{\gamma + \sqrt{\gamma^2 + 4\theta}}{2\theta}$ where $\theta$ is the constant of semi-convexity of $V$ and $\gamma$ is the friction coefficient.

**Lemma 2.7** *Let $(q, p) \in \mathbb{R}^{2d}$ and $\delta < \delta_0 := \frac{\gamma + \sqrt{\gamma^2 + 4\theta}}{2\theta}$. Under the assumption of semi-convexity of $V$ (**B-1**) and the dissipativity inequality (2.3), there exists a unique $z \in \mathbb{R}^d$ such that $z = q + \frac{\delta}{1 + \gamma \delta}(p - \delta \partial_q V(z))$.*

*Proof* Let $(q,p) \in \mathbb{R}^{2d}$ and $\delta < \delta_0$. Let $P : \mathbb{R}^d \to \mathbb{R}^d$ defined for $z \in \mathbb{R}^d$ by $P(z) = -z + q + \frac{\delta}{1+\gamma\delta}(p - \delta\partial_q V(z))$. We have that $P \in C^\infty$. Using dissipativity inequality (2.3), we get for all $z \in \mathbb{R}^d$

$$\langle P(z), z \rangle = -|z|^2 + \langle q, z \rangle - \frac{\delta^2}{1+\delta\gamma}\langle \partial_q V(z), z \rangle + \frac{\delta}{1+\delta\gamma}\langle p, z \rangle$$

$$\leq -|z|^2(1 + \frac{\beta_1 \delta^2}{1+\delta\gamma}) + \frac{\delta}{1+\delta\gamma}\langle p, z \rangle + \kappa\frac{\delta^2}{1+\delta\gamma} + \langle q, z \rangle.$$

Thus, we have for $|z|^2$ large enough than $\langle P(z), z \rangle < 0$. Then, using a Corollary of Brower fixed-point theorem (see for instance [17]), we have that there exists $z \in \mathbb{R}^d$ such that $P(z) = 0$. Therefore we have shown the existence of $z \in \mathbb{R}^d$ such that $z = q + \frac{\delta}{1+\gamma\delta}(p - \delta\partial_q V(z))$.

Let us show the uniqueness.
Let $z_1 \in \mathbb{R}^d$ and $z_2 \in \mathbb{R}^d$ such that $P(z_1) = P(z_2) = 0$ and $z_1 \neq z_2$. We have

$$z_1 - z_2 = -\frac{\delta^2}{1+\delta\gamma}(\partial_q V(z_1) - \partial_q V(z_2)).$$

Using assumption of semi-convexity (**B-1**), we get

$$|z_1 - z_2|^2 = -\frac{\delta^2}{1+\delta\gamma}\int_0^1 \partial_q^2 V(z_1 + t(z_2 - z_1))(z_1 - z_2, z_1 - z_2)dt \leq \frac{\delta^2 \theta}{1+\delta\gamma}|z_1 - z_2|^2.$$

Since $\delta < \frac{\gamma + \sqrt{\gamma^2 + 4\theta}}{2\theta}$, we have $\frac{\delta^2 \theta}{1+\delta\gamma} < 1$ and $z_1 = z_2$.

**Remark 2.8** *The condition $\delta < \delta_0 := \frac{\gamma + \sqrt{\gamma^2 + 4\theta}}{2\theta}$ is useful only for the uniqueness. We have the existence for all $\delta$. Moreover, in the case where V is convex, the inequality on the second derivative is true for all $\theta \geq 0$. Hence, we can show the uniqueness for all $\delta > 0$.*
*An other proof of the fact that the implicit split-step scheme is well defined can be found in [18].*

Moreover for $\delta$ small enough, we have that these two schemes have moments of all order. More exactly, we have:

**Proposition 2.9** *Let $k \in \mathbb{N}^*$ and $\delta < \delta_0 := \frac{\gamma\beta}{4\theta}$. Let $(q,p) \in \mathbb{R}^{2d}$ such that $q_0 = q(0)$ and $p_0 = p(0)$. Under assumptions **B**, the implicit split-step scheme (2.11) and the implicit Euler scheme (2.12) satisfy: There exist positive numbers $C_k$ and $\ell_k$ such that for all $n \in \mathbb{N}$*

$$\mathbb{E}(|q_n|^{2k} + |p_n|^{2k}) < C_k(1 + |q(0)|^{\ell_k} + |p(0)|^{\ell_k}). \tag{2.13}$$

**Remark 2.10** *Since $\frac{\gamma\beta}{4\theta} < \frac{\gamma + \sqrt{\gamma^2 + 4\theta}}{2\theta}$, the schemes considered are well defined. If we assume that $\delta_0 < 1$ then we can choose $C_k$ independent of $\delta_0$.*
*We will prove this result only for the implicit Euler scheme. The proof for the implicit split-step scheme is similar. A proof for the moment of order 2 of the implicit split-step scheme can be found in [18]. The two proofs use the same ideas.*

To prove Proposition 2.9 for the implicit Euler scheme (2.12), we need the following three Lemmas:

**Lemma 2.11** *Let $V \in C^\infty(\mathbb{R}^d)$. Let us assume that V is semi-convex, then V verifies for any $q \in \mathbb{R}^d$ and $p \in \mathbb{R}^d$*

$$V(q) - V(p) \leq \langle \partial_q V(q), q - p \rangle + \frac{\theta}{2}|q - p|^2, \tag{2.14}$$

*where $\theta$ is the constant of semi-convexity of V.*

*Proof* Let $q \in \mathbb{R}^d$ and $p \in \mathbb{R}^d$. Using the Taylor expansion on $V$ and the semi-convexity assumption on $V$ (**B-1**), we have

$$V(p) - V(q) = \langle \partial_q V(q), p - q \rangle + \int_0^1 (1-s) \partial_q^2 V(q + s(p-q))(p-q, p-q) ds$$
$$\geq \langle \partial_q V(q), p - q \rangle - \frac{\theta}{2} |q - p|^2.$$

Finally, we get

$$V(q) - V(p) \leq \langle \partial_q V(q), q - p \rangle + \frac{\theta}{2} |q - p|^2.$$

**Lemma 2.12** *Let $0 < \delta < \frac{\gamma \beta}{4\theta}$ and $\Gamma_\delta \in C^\infty(\mathbb{R}^{2d})$ be defined for any $(q,p) \in \mathbb{R}^{2d}$ by*

$$\Gamma_\delta(q,p) = \Gamma(q,p) + \frac{\gamma \delta}{4} |p|^2$$
$$= \frac{1}{2} |p|^2 + V(q) + \frac{\gamma}{2} \langle q, p \rangle + \frac{\gamma^2}{4} |q|^2 + 1 + \frac{\gamma \delta}{4} |p|^2, \tag{2.15}$$

*then $\Gamma_\delta$ verifies the following properties: for any $(q,p) \in \mathbb{R}^{2d}$*

$$\Gamma_\delta(q,p) \geq \frac{1}{8} |p|^2 \tag{2.16}$$

*and for any $n \in \mathbb{N}$*

$$(1 + \gamma(1-\varepsilon)\delta\beta)\Gamma_\delta(q_{n+1}, p_{n+1}) \leq \Gamma_\delta(q_n, p_{n+1} + \delta\gamma p_{n+1} + \delta \partial_q V(q_{n+1})) + \gamma\delta(\kappa + \beta), \tag{2.17}$$

*where $(q_n, p_n)$ is the implicit Euler scheme defined by (2.12) and $0 < \varepsilon < 1$ is such that $\delta < \frac{\varepsilon \gamma \beta}{4\theta}$.*

*Proof* Inequality (2.16) is a corollary of property (2.6) on $\Gamma$.
Using the definition of the implicit Euler scheme (2.12), we have

$$\Gamma_\delta(q_{n+1}, p_{n+1}) = \frac{1}{2} |p_{n+1}|^2 + V(q_n + \delta p_{n+1}) + \frac{\gamma}{2} \langle q_n, p_{n+1} \rangle + \frac{\gamma \delta}{2} |p_{n+1}|^2$$
$$+ \frac{\gamma^2}{4} \left( |q_n|^2 + 2\delta \langle q_n, p_{n+1} \rangle + \delta^2 |p_{n+1}|^2 \right) + 1 + \frac{\gamma \delta}{4} |p_{n+1}|^2.$$

Using inequality (2.14) on $V$, we get

$$\Gamma_\delta(q_{n+1}, p_{n+1}) \leq \Gamma_\delta(q_n, p_{n+1}) + \delta \langle \partial_q V(q_{n+1}), p_{n+1} \rangle + \left( \frac{\delta\gamma}{2} + \frac{\theta\delta^2}{2} + \frac{\gamma^2 \delta^2}{4} \right) |p_{n+1}|^2$$
$$+ \frac{\gamma^2 \delta}{2} \langle q_n, p_{n+1} \rangle.$$

Since $q_n = q_{n+1} - \delta p_{n+1}$ and, for any $(q,p) \in \mathbb{R}^{2d}$,

$$|p|^2 - |q|^2 = \langle p-q, p+q \rangle = 2\langle p-q, p \rangle - |p-q|^2 \leq 2\langle p-q, p \rangle,$$

we get

$$\begin{aligned}\Gamma_\delta(q_n,p_{n+1}) =& \Gamma_\delta(q_n,p_{n+1}+\delta\gamma p_{n+1}+\delta\partial_q V(q_{n+1}))\\ &+\frac{1}{2}\big(|p_{n+1}|^2-|(1+\delta\gamma)p_{n+1}+\delta\partial_q V(q_{n+1})|^2\big)\\ &-\frac{\gamma\delta}{2}\langle q_{n+1},\partial_q V(q_{n+1})\rangle-\frac{\gamma^2\delta}{2}\langle q_{n+1},p_{n+1}\rangle\\ &+\frac{\gamma\delta^2}{2}\langle p_{n+1},\partial_q V(q_{n+1})+\gamma p_{n+1}\rangle\\ &+\frac{\gamma\delta}{4}\big(|p_{n+1}|^2-|(1+\delta\gamma)p_{n+1}+\delta\partial_q V(q_{n+1})|^2\big)\\ \leq & \Gamma_\delta(q_n,p_{n+1}+\delta\gamma p_{n+1}+\delta\partial_q V(q_{n+1}))-\gamma\delta|p_{n+1}|^2-\delta\langle\partial_q V(q_{n+1}),p_{n+1}\rangle\\ &-\frac{\gamma\delta}{2}\langle q_{n+1},\partial_q V(q_{n+1})\rangle-\frac{\gamma^2\delta}{2}\langle q_{n+1},p_{n+1}\rangle.\end{aligned}$$

Then, we get

$$\begin{aligned}\Gamma_\delta(q_{n+1},p_{n+1}) \leq & \Gamma_\delta(q_n,p_{n+1}+\delta\gamma p_{n+1}+\delta\partial_q V(q_{n+1}))+\big(\frac{\theta\delta^2}{2}-\frac{\gamma^2\delta^2}{4}\big)|p_{n+1}|^2\\ &-\frac{\delta\gamma}{2}\big(\langle q_{n+1},\partial_q V(q_{n+1})\rangle+|p_{n+1}|^2\big).\end{aligned}$$

Using (2.8) and $\beta<1$, we have

$$\begin{aligned}\Gamma_\delta(q_n,p_{n+1}) \leq & \Gamma_\delta(q_n,p_{n+1}+\delta\gamma p_{n+1}+\delta\partial_q V(q_{n+1}))-\gamma\delta\beta\Gamma_\delta(q_{n+1},p_{n+1})\\ &+\big(\frac{\theta\delta^2}{2}+\frac{\gamma^2\delta^2}{4}(\beta-1)\big)|p_{n+1}|^2+(\kappa+\beta)\gamma\delta\\ \leq & \Gamma_\delta(q_n,p_{n+1}+\delta\gamma p_{n+1}+\delta\partial_q V(q_{n+1}))-\gamma\delta\beta\Gamma_\delta(q_{n+1},p_{n+1})\\ &+\frac{\theta\delta^2}{2}|p_{n+1}|^2+(\kappa+\beta)\gamma\delta.\end{aligned}$$

Using property (2.16), we get

$$\begin{aligned}\Gamma_\delta(q_{n+1},p_{n+1}) \leq & \Gamma_\delta(q_n,p_{n+1}+\delta\gamma p_{n+1}+\delta\partial_q V(q_{n+1}))+\big(\frac{\theta\delta^2}{2}-\beta\frac{\gamma\delta\varepsilon}{8}\big)|p_{n+1}|^2\\ &+(\kappa+\beta)\gamma\delta-\gamma\delta\beta(1-\varepsilon)\Gamma_\delta(q_{n+1},p_{n+1}).\end{aligned}$$

Since $\delta\leq\frac{\varepsilon\gamma\beta}{4\theta}$ it follows that

$$(1+\gamma(1-\varepsilon)\delta\beta)\Gamma_\delta(q_{n+1},p_{n+1}) \leq \Gamma_\delta(q_n,p_{n+1}+\delta\gamma p_{n+1}+\delta\partial_q V(q_{n+1}))+\gamma\delta(\kappa+\beta).$$

**Lemma 2.13** *Let $\delta<\delta_0:=\frac{\gamma+\sqrt{\gamma^2+4\theta}}{2\theta}$ be fixed. Let the processes $P_n(.)$ be defined for $t\in[n\delta,(n+1)\delta]$ and $n\in\mathbb{N}^*$ by*

$$P_n(t)=p_n+\sigma\big(W(t)-W(n\delta)\big).$$

*Then, there exists $C(\delta_0)>0$ such that for all $n\in\mathbb{N}$ and $t\in[n\delta,(n+1)\delta]$*

$$\mathbb{E}\big(\Gamma_\delta^\ell(q_n,P_n(t))\big)\leq\mathbb{E}\big(\Gamma_\delta^\ell(q_n,p_n)\big)+C(\delta_0)\sum_{i=0}^{\ell-1}(t-n\delta)^{\ell-i}\mathbb{E}\big(\Gamma_\delta^i(q_n,p_n)\big). \tag{2.18}$$

*Proof* We will show (2.18) by induction on $\ell$. Let $n \in \mathbb{N}^*$ be fixed. Using the definition of $\Gamma_\delta$, independence of $q_n$ with $W(t) - W(n\delta)$ for $t \geq n\delta$ and properties of $W$, we easily have

$$E\left(\Gamma_\delta(q_n, p_n + \sigma(W(t) - W(n\delta)))\right) \leq \mathbb{E}\left(\Gamma_\delta(q_n, p_n)\right) + (t - n\delta)\frac{\sigma^2}{2}(1 + \frac{\delta_0 \gamma}{2}).$$

Let $\ell \in \mathbb{N}^*$. Let us assume that (2.18) is true for all $j \leq \ell - 1$ and let us show it for $\ell$. Let $n \in \mathbb{N}^*$ be fixed. Using the Itô formula on $\Gamma_\delta(q_n, P_n(t))$, we get for $t \in [n\delta, (n+1)\delta]$

$$\Gamma_\delta^\ell(q_n, P_n(t)) = \Gamma_\delta^\ell(q_n, p_n) + \sigma \ell \int_{n\delta}^t \Gamma_\delta^{\ell-1}(q_n, P_n(s)) \langle (1 + \frac{\gamma\delta}{2}) P_n(s) + \frac{\gamma}{2} q_n, dW(s) \rangle$$
$$+ \frac{1}{2}\ell(1 + \frac{\delta\gamma}{2}) \int_{n\delta}^t \Gamma_\delta^{\ell-1}(q_n, P_n(s)) ds$$
$$+ \frac{\ell(\ell-1)}{2} \int_{n\delta}^t |(1 + \frac{\gamma\delta}{2}) P_n(s) + \frac{\gamma}{2} q_n|^2 \Gamma_\delta^{\ell-2}(q_n, P_n(s)) ds.$$

We take the expectation. The second term in the right hand side vanishes because it is a square integrable martingale. Using for $(q, p) \in \mathbb{R}^{2d}$

$$\left|(1 + \frac{\gamma\delta}{2})p + \frac{\gamma}{2}q\right|^2 \leq 16(1 + \frac{\gamma\delta_0}{2})^2 \Gamma_\delta(q, p),$$

we get, for $t \in [n\delta, (n+1)\delta]$,

$$\mathbb{E}\left(\Gamma_\delta^\ell(q_n, P_n(t))\right) \leq \mathbb{E}\left(\Gamma_\delta^\ell(q_n, p_n)\right) + C_\ell(\delta_0) \int_{n\delta}^t \mathbb{E}\left(\Gamma_\delta^{\ell-1}(q_n, P_n(s))\right) ds.$$

Using the induction hypothesis, we get (2.18) for $\ell$.

*Proof (Proof of Proposition 2.9 in the case of the implicit Euler scheme (2.12).)* First, we rewrite the implicit Euler scheme (2.12) like this

$$\begin{aligned} q_n &= q_{n+1} - \delta p_{n+1} \\ p_n + \sqrt{\delta}\sigma \eta_n &= p_{n+1} + \delta \partial_q V(q_{n+1}) + \gamma \delta p_{n+1}. \end{aligned} \quad (2.19)$$

Let $\ell \in \mathbb{N}^*$. Using (2.18) for $t = (n+1)\delta$, we see that

$$\mathbb{E}\left(\Gamma_\delta^\ell(q_n, p_n + \sigma\sqrt{\delta}\eta_n)\right) \leq \mathbb{E}\left(\Gamma_\delta^\ell(q_n, p_n)\right) + C \sum_{i=0}^{\ell-1} \delta^{\ell-i} \mathbb{E}\left(\Gamma_\delta^i(q_n, p_n)\right)$$
$$\leq \mathbb{E}\left(\Gamma_\delta^\ell(q_n, p_n)\right) + C\delta \sum_{i=0}^{\ell-1} \mathbb{E}\left(\Gamma_\delta^i(q_n, p_n)\right),$$

where we can choose $C$ independent of $\delta$. Using the fact that for any $x \in \mathbb{R}$, $i \in \mathbb{N}^*$ such that $i < \ell$ and $\varepsilon_1 > 0$

$$|x|^i \leq \varepsilon_1 |x|^\ell + C_{\varepsilon_1}, \quad (2.20)$$

we get

$$\mathbb{E}\left(\Gamma_\delta^\ell(q_n, p_n + \sigma\sqrt{\delta}\eta_n)\right) \leq (1 + \varepsilon_1 \delta C) \mathbb{E}\left(\Gamma_\delta^\ell(q_n, p_n)\right) + CC_{\varepsilon_1}\delta \quad (2.21)$$

Moreover, using (2.20) and (2.17), we get that for $\varepsilon_2 > 0$

$$(1 + \gamma(1 - \varepsilon)\delta\beta)^\ell \Gamma_\delta^\ell(q_{n+1}, p_{n+1}) \leq (1 + \delta\varepsilon_2) \Gamma_\delta^\ell(q_n, p_{n+1} + \delta\gamma p_{n+1} + \delta\partial_q V(q_{n+1})) + \delta C_{\varepsilon_2}, \quad (2.22)$$

where $0 < \varepsilon < 1$ is defined in (2.17). Then, using (2.21), (2.19) and (2.22), we get

$$(1+\gamma(1-\varepsilon)\delta\beta)^\ell \mathbb{E}\big(\Gamma_\delta^\ell(q_{n+1},p_{n+1})\big) \leq (1+\delta\varepsilon_2)(1+\varepsilon_1\delta C)\mathbb{E}\big(\Gamma_\delta^\ell(q_n,p_n)\big) + \delta C_{\varepsilon_1,\varepsilon_2}.$$

Choosing $\varepsilon_1 = \frac{\gamma(1-\varepsilon)\beta}{2C}$ and $\varepsilon_2 = \lambda\gamma(1-\varepsilon)\beta$ such that $\lambda \leq \frac{1}{2(1+1/2(1-\varepsilon)\gamma\beta\delta_0)}$, we get, $(1+\varepsilon_1 C)(1+\varepsilon_2) \leq 1+\gamma(1-\varepsilon)\beta$. Then, by induction on $n$, we get

$$\mathbb{E}\big(\Gamma_\delta^\ell(q_n,p_n)\big) \leq \Gamma_\delta^\ell(q,p) + C\delta \sum_{i=1}^n \frac{1}{(1+\gamma(1-\varepsilon)\delta\beta)^{(\ell-1)i}}$$
$$\leq \Gamma_\delta^\ell(q,p) + C.$$

To conclude, we use the fact the $\Gamma_\delta$ has polynomial growth and for any $(q,p) \in \mathbb{R}^{2d}$

$$\frac{1}{8}|p|^2 + \frac{\gamma^2}{12}|q|^2 \leq c\Gamma_\delta(q,p).$$

2.3 Preliminary results

For $(x,y) \in \mathbb{R}^{2d}$, we denote by $\big(q_x(t), p_y(t)\big)_{t\geq 0}$ a process which verifies (2.4) and has for initial data $q_x(0) = x$ and $p_y(0) = y$. From now on, $(P_t)_{t\geq 0}$ is the transition semigroup associated with the Markov process $\big(q_x(t), p_y(t)\big)_{t\geq 0}$.

We recall that we denote by $L$ the Kolmogorov generator associated with the stochastic equation (2.4) defined for all $\phi \in C^\infty(\mathbb{R}^{2d})$ and $(q,p) \in \mathbb{R}^{2d}$ by

$$L\phi(q,p) := \langle p, \partial_q \phi(q,p)\rangle - \langle \partial_q V(q) + \gamma p, \partial_p \phi(q,p)\rangle + \frac{\sigma^2}{2}\sum_{i=1}^d \frac{\partial^2}{\partial p_i \partial p_i}\phi(q,p). \quad (2.23)$$

Moreover its formal adjoint $L^\top$ in $\mathbb{R}^{2d}$ is defined for all $\phi \in C^\infty(\mathbb{R}^{2d})$ and $(q,p) \in \mathbb{R}^{2d}$ by

$$L^\top \phi(q,p) = -\langle p, \partial_q \phi(q,p)\rangle + \langle \partial_q V(q), \partial_p \phi(q,p)\rangle + \gamma\langle p, \partial_p \phi(q,p)\rangle$$
$$+ \frac{\sigma^2}{2}\sum_{i=1}^d \frac{\partial^2}{\partial p_i \partial p_i}\phi(q,p) + d\gamma\phi(q,p). \quad (2.24)$$

The following equality will be useful: For any functions $\phi \in C^\infty(\mathbb{R}^{2d})$ and $\psi \in C^\infty(\mathbb{R}^{2d})$, we have

$$L(\phi\psi) = \psi L\phi + \phi L\psi + \sigma^2 \langle \partial_p \phi, \partial_p \psi\rangle. \quad (2.25)$$

We define the measure $d\rho := \rho(q,p)dqdp$ where for all $(q,p) \in \mathbb{R}^{2d}$, $\rho(q,p) = \frac{1}{Z}\exp(-\frac{2\gamma}{\sigma^2}H(q,p))$ and $Z = \int_{\mathbb{R}^{2d}} \exp(-\frac{2\gamma}{\sigma^2}H(q,p))dqdp$. A consequence of assumption **B-4** is that for all $k \in \mathbb{N}$ and $j \in \mathbb{N}$

$$\int_{\mathbb{R}^{2d}} |q|^{2k}|p|^{2j} e^{-\frac{2\gamma}{\sigma^2}H(q,p)}dqdp < \infty.$$

It is easy to verify that the measure $d\rho$ is invariant by $P_t$: $L^\top \rho = 0$.

In all this article, we use the following spaces:

$$L^2(\rho) = \{f : \mathbb{R}^{2d} \to \mathbb{R}; \int |f|^2 d\rho < \infty\},$$

and for $j \in \mathbb{N}$ and $\ell \in \mathbb{N}$

$$\mathscr{C}_j^\ell(\mathbb{R}^{2d}) := \Big\{ f \in C^\ell(\mathbb{R}^{2d}) \text{ such that for all } \mathbf{k} = (k_1,...,k_{2d}) \in \mathbb{N}^d, \ |\mathbf{k}| \leq j,$$
$$\sup_{(q,p) \in \mathbb{R}^{2d}} |\partial^\mathbf{k} f(q,p)| \Gamma^{-\ell}(q,p) < \infty \Big\}.$$

Let $\ell$ and $k$ be two integers. Let a function $f \in \mathscr{C}_k^\ell(\mathbb{R}^{2d})$, we define the following norm

$$\| f \|_{\ell,k} := \sup_{\substack{\mathbf{j}:=(j_1,...j_{2d}) \\ |\mathbf{j}| \leq k}} \sup_{x \in \mathbb{R}^{2d}} \left| \frac{\partial^\mathbf{j} f(x)}{\Gamma^\ell(x)} \right|.$$

We also define the semi-norm

$$|f|_{\ell,k} := \sup_{\substack{\mathbf{j}:=(j_1,...j_{2d}) \\ 1 \leq |\mathbf{j}| \leq k}} \sup_{x \in \mathbb{R}^{2d}} \left| \frac{\partial^\mathbf{j} f(x)}{\Gamma^\ell(x)} \right|.$$

We consider a function $\phi \in \mathscr{C}_{pol}^\infty(\mathbb{R}^{2d})$ and we set for all $(q,p) \in \mathbb{R}^{2d}$ and $t \geq 0$:

$$u(t,q,p) = \mathbb{E}(\phi(q_q(t), p_p(t))). \tag{2.26}$$

We have that $u$ is a $C^\infty$ function on $\mathbb{R}^+ \times \mathbb{R}^{2d}$. Moreover we have that $u$ is the unique solution of the Kolmogorov equation

$$\frac{du}{dt}(t,q,p) = Lu(q,p), \ (q,p) \in \mathbb{R}^{2d}, \ t > 0, \quad u(0,q,p) = \phi(q,p), \ (q,p) \in \mathbb{R}^{2d}. \tag{2.27}$$

In the following, we write: $u(t) = P_t \phi$. Note that we use the standard identification $u(t) = u(t,.)$. Moreover, we have the following result on the regularity of $u$ (see [26]):

**Lemma 2.14** *Let $\phi \in \mathscr{C}_{pol}^\infty(\mathbb{R}^{2d})$. Let $u$ defined by (2.26). For any integer, there exists an integer $s$ such that for all $T \geq 0$ there exists a strictly positive real number $C_m(T)$ such that for $t \in [0,T]$ and $(q,p) \in \mathbb{R}^{2d}$*

$$\left| D^m u(t,q,p) \right| \leq C_m(T)(1 + |q|^s + |p|^s),$$

*where $D^m u(t)$ denotes the vector of all the spatial derivatives of $u(t)$ of order $m$.*

We wish to investigate the approximation properties of the implicit Euler scheme and the implicit split-step scheme. We need results on the long time behavior of the law of the solution (2.27). The two necessary properties are the following :

**Proposition 2.15** *Let $\phi \in \mathscr{C}_{pol}^\infty(\mathbb{R}^{2d})$ such that $\int \phi d\rho = 0$. Let $u$ the unique solution of (2.27). There exists a strictly positive real number $\lambda_0$ such that, for $m \in \mathbb{N}$ and $\lambda < \lambda_0$, there exist a positive real number $C$ and integers $r_{m+d+1}$ and $\ell_m \geq 2r_{m+d+1}$ such that $\phi \in \mathscr{C}_{r_{m+d+1}}^{m+d+1}(\mathbb{R}^{2d})$ and for all $t \geq 0$*

$$\| u(t) \|_{m,\ell_m} \leq C \exp(-\lambda t) \| \phi \|_{d+1+m, r_{m+d+1}}.$$

**Lemma 2.16** *We denote by $L^*$ the formal adjoint of $L$ in $L^2(\rho)$. Let $g \in \mathscr{C}_{pol}^\infty(\mathbb{R}^{2d})$ such that $\int g d\rho = 0$. Then there exists a unique function $\mu \in \mathscr{C}_{pol}^\infty(\mathbb{R}^{2d})$ such that*

$$L^* \mu = g \quad \text{and} \quad \int \mu d\rho = 0.$$

The proof of this two results can be found in the appendix.

Our main result can be stated as follows:

**Theorem 2.17** *Let N be fixed. Let $\delta_0 := \min(\frac{1}{\gamma}, \frac{\gamma\beta}{4\theta})$. Let $(q_k, p_k)$ be the discrete process defined by the implicit Euler scheme (2.12) or by the implicit split-step scheme (2.11), then for all $\delta < \delta_0$, there exists a modified function $\mu^{(N)}$ defined for all $(q,p) \in \mathbb{R}^{2d}$ by*

$$\mu^{(N)}(q,p) = 1 + \sum_{n=1}^{N} \delta^n \mu_n(q,p)$$

*such that $\mu^{(N)} \in \mathscr{C}_{pol}^{\infty}(\mathbb{R}^{2d})$ and*

$$\int_{\mathbb{R}^{2d}} \mu^{(N)}(q,p)\rho(q,p)dqdp = 1.$$

*For all function $\phi \in \mathscr{C}_{pol}^{\infty}(\mathbb{R}^{2d}) \cap \mathscr{C}_{\ell_N}^{\alpha_N}(\mathbb{R}^{2d})$ where $\alpha_N = 6N+2+(d+1)(N+1)$ and $\lambda < \lambda_0$, where $\lambda_0$ is defined in Proposition 2.15, there exist a positive real number $C_N$, an integer $k_N \geq \ell_N$ and a positive polynomial function $P_N$ satisfying the following : For all $k \in \mathbb{N}$,*

$$\| \mathbb{E}\phi(q_k, p_k) - \int \phi \mu^{(N)} d\rho \|_{0,k_N} \leq C_N \left(e^{-\lambda t_k} + \delta^N\right) \| \phi - \langle \phi \rangle \|_{\alpha_N, \ell_N},$$

*where $t_k = k\delta$ and $\langle \phi \rangle = \int \phi d\rho$.*

This result can be viewed as a discrete version of Proposition 2.15 in the case $m = 0$. We have, for $(q_k, p_k)$ the discrete process defined by the implicit Euler scheme (2.12) or the implicit split-step scheme (2.11), that $\mathbb{E}\phi(q_k, p_k)$, which is an approximation of $u$, has the same property as $u$ : $\mathbb{E}\phi(q_k, p_k)$ converge exponentially fast to a constant in $\mathscr{C}_{k_N}^{0}(\mathbb{R}^{2d})$ up to an error $\delta^N C_N$. At $k$ fixed, we can optimize this error with a good choice of $N$.

Our result can be compared with [4,12,19,24,25]. As in [4,12,19], the only assumption made on $\delta$ is that $\delta < \delta_0 := \min(\frac{1}{\gamma}, \frac{\gamma\beta}{4\theta})$. We also recover an expansion of the invariant measure as in [27].

Our result is similar to the result in the case of SDE on the torus describes in [4] and to the result for the overdamped Langevin equation describes in [12].

The constant $C_N$ appearing in the estimate depends of $N$, the constant of semi-convexity $\theta$, the parameters of the equation $\gamma$ and $\sigma$ and the polynomial growth of $V$ and all its derivatives.

## 3 Asymptotic expansion of the weak error

We have the formal expansion for small $t$ and $(q,p) \in \mathbb{R}^{2d}$:

$$u(t,q,p) = \phi(q,p) + tL\phi(q,p) + \frac{t^2}{2}L^2\phi(q,p) + ... + \frac{t^n}{n!}L^n\phi(q,p) + ...$$

This is just obtained by Taylor expansion in time.

Since the solution $u(t)$ of the Kolmogorov equation is in $\mathscr{C}_{pol}^{\infty}(\mathbb{R}^{2d})$, the above formal expansion can be justified in $\mathscr{C}_{pol}^{\infty}(\mathbb{R}^{2d})$. Indeed, we have the following easy result whose proof is left to the reader.

**Proposition 3.1** *Let $\phi \in \mathscr{C}_{pol}^{\infty}(\mathbb{R}^{2d})$ and $\delta_1 > 0$ be fixed. Then, for all $N$, there exist constants $C(N)$ and $\ell$ such that for all $\delta < \delta_1$,*

$$|u(\delta,q,p) - \sum_{n=0}^{N} \frac{\delta^n}{n!} L^n(q,p)\phi(q,p)| \leq C(N)\delta^{N+1}(1+|q|^{\ell}+|p|^{\ell}) \| \phi \|_{2N+2, r_{2N+2}},$$

*where $r_{2N+2}$ is defined such that $\phi \in \mathscr{C}_{r_{2N+2}}^{2N+2}(\mathbb{R}^{2d})$.*

We now examine in detail the first time step and its approximation properties in terms of law. By Markov property, it is sufficient to then obtain information at all steps. We want to have an expansion similar to the last Proposition for the processes defined by the implicit Euler scheme (2.12) or the implicit split-step scheme (2.11).

**Proposition 3.2** *Let $\delta_0 = \min(\frac{1}{\gamma}, \frac{\gamma\beta}{4\theta})$. Let $\phi \in \mathscr{C}^\infty_{pol}(\mathbb{R}^{2d})$. For any $N \in \mathbb{N}$, there exists an integer $r_N$ such that $\phi \in \mathscr{C}^N_{r_N}(\mathbb{R}^{2d})$. For all $n \geq 1$, there exist operators $A_n$ of order $2n$ with coefficients $\mathscr{C}^\infty_{pol}(\mathbb{R}^{2d})$ which depend of the scheme chosen ((2.11) or (2.12)), such that for all integer $N \geq 1$ there exist a constant $C_N$ and an integer $\ell_N$ such that for $0 < \delta < \delta_0$ and $(q,p) := (q_0, p_0)$,*

$$|\mathbb{E}\phi(q_1, p_1) - \sum_{n=0}^{N} \delta^n A_n(q,p)\phi(q,p)| \leq C_N \delta^{N+1}(1 + |q|^{\ell_N} + |p|^{\ell_N})|\phi|_{2N+2, r_{2N+2}}. \tag{3.1}$$

*Moreover, we have $A_0 = I$ and $A_1 = L$.*

**Remark 3.3** *This result is similar of the asymptotic expansion of the weak error describes in [4], but the proof is different. Indeed, we can not use Itô's lemma because the schemes considered here are implicit.*

We consider the implicit split-step scheme (2.11). Let $0 < \delta < \delta_0$. Let us recall that we have

$$\begin{aligned} q_1 &= q_0 + \delta p_0^* \\ p_0^* &= p_0 - \delta\gamma p_0^* - \delta\partial_q V(q_1) \\ p_1 &= p_0^* + \sqrt{\delta}\sigma\eta_0. \end{aligned}$$

Before proving Proposition 3.2, we need an asymptotic expansion for $q_1$, $p_1$ and $p_0^*$.
We define the function $\Psi_\delta$ which associate to $(q,p)$ the solution $z$ of $(1+\gamma\delta)z = (1+\gamma\delta)q + \delta p - \delta^2\partial_q V(z)$. The function $\Psi_\delta$ is well defined (see Lemma 2.7) and we have $q_1 = \Psi_\delta(q_0, p_0)$. Moreover, we have that $(\delta, q, p) \mapsto \Psi_\delta(q, p)$ is $C^\infty$ on $]0, \frac{\gamma+\sqrt{\gamma^2+4\theta}}{2\theta}[\times\mathbb{R}^{2d}$: Let $\Omega_1 = ]0, \frac{\gamma+\sqrt{\gamma^2+4\theta}}{2\theta}[\times\mathbb{R}^{3d}$ and the function $f \in C^\infty$ defined on $\Omega_1$ by

$$f(\delta, q, p, z) = -(1+\delta\gamma)z + (1+\delta\gamma)q + \delta p - \delta^2\partial_q V(z).$$

Using the semi-convexity of $V$, we have that, for all $(\delta, q, p, z) \in \Omega_1$, $\partial_z f(\delta, q, p, z)$ is invertible. Using implicit function Theorem, we have that the function defined by $(\delta, q, p) \mapsto \Psi_\delta(q, p) = z$ is $C^\infty$ on a neighborhood of each point of $]0, \frac{\gamma+\sqrt{\gamma^2+4\theta}}{2\theta}[\times\mathbb{R}^{2d}$.

We have the following result:

**Lemma 3.4** *Let $\delta_0 = \min(\frac{1}{\gamma}, \frac{\gamma\beta}{4\theta})$. Let $(q,p) \in \mathbb{R}^{2d}$ such that $q_0 = q$ and $p_0 = p$. We have for $0 < \delta < \delta_0$ and $N \in \mathbb{N}$,*

$$q_1 = \Psi_\delta(q,p) = q + \sum_{k=1}^{N} \delta^k d_k(q,p) + \delta^{N+1} R_{N+1}(q, p, \delta), \tag{3.2}$$

*where, for all $k \geq 0$, $d_k \in \mathscr{C}^\infty_{pol}(\mathbb{R}^{2d})$ is defined for all $(x,y) \in \mathbb{R}^{2d}$ by*

$$\begin{aligned} d_1(x,y) &= y, \text{ and for } k \geq 2 \\ d_k(x,y) &= (-1)^{k-1}\gamma^{k-2}\big(\gamma y + \partial_q V(x)\big) \\ &\quad + \sum_{j=2}^{k-1}(-1)^{j-1}\gamma^{j-2}\sum_{n=1}^{k-j}\frac{1}{n!}\sum_{\substack{k_1+\ldots+k_n=k-j-n,\\ 0\leq k_i \leq N}} \partial_q^{n+1} V(x)(d_{k_1+1}(x,y), \ldots, d_{k_n+1}(x,y)). \end{aligned}$$

*and $R_{N+1}$ verifies: There exist $C > 0$ and $\ell_N \in \mathbb{N}$ such that for any $(x, y) \in \mathbb{R}^{2d}$ and $\delta < \delta_0$*

$$|R_{N+1}(x, y, \delta)| \leq C(1 + |x|^{\ell_N} + |y|^{\ell_N}). \tag{3.3}$$

*Proof* Let $0 < \delta < \delta_0$. We have previously shown that $(\delta, q, p) \mapsto \Psi_\delta(q, p)$ is $C^\infty$ on $]0, \frac{\gamma + \sqrt{\gamma^2 + 4\theta}}{2\theta}[ \times \mathbb{R}^{2d}$. Then, for $(q, p) \in \mathbb{R}^{2d}$ fixed, we have that $d_k(q, p)$ is the $k^{th}$ term of the Taylor expansion of $\delta \mapsto \Psi_\delta(q, p)$ and we can write (3.2). We now search an expression for $d_k$.

Let $(q, p) \in \mathbb{R}^{2d}$ such that $q_0 = q$ and $p_0 = p$. Let $\ell \in \mathbb{N}$. We use the temporary notation, for all $0 \leq k \leq \ell - 1$, $g_{k,\ell} = d_{k+1}$ and $g_{\ell,\ell} = R_{\ell+1}$. Hence, we have

$$q_1 = q + \sum_{k=1}^{\ell} \delta^k d_k(q, p) + \delta^{\ell+1} R_{\ell+1}(q, p, \delta)$$
$$= q + \delta z_{1,\ell},$$

where $z_{1,\ell} = \sum_{k=0}^{\ell} g_{k,\ell}(q, p)$. Using Taylor expansion, we obtain

$$\partial_q V(q_1) = \partial_q V(q + \delta z_{1,l}) = \partial_q V(q) + \sum_{n=1}^{\ell} \frac{1}{n!} \partial_q^{n+1} V(q)(\delta z_{1,l}, ..., \delta z_{1,l}) + \delta^{\ell+1} \theta_\ell(q, p)$$

$$= \partial_q V(q) + \sum_{n=1}^{\ell} \frac{1}{n!} \delta^n \partial_q^{n+1} V(q) (\sum_{k=0}^{\ell} \delta^k g_{k,\ell}(q, p), ..., \sum_{k=0}^{\ell} \delta^k g_{k,\ell}(q, p)) + \delta^{\ell+1} \theta_\ell(q, p)$$

$$= \partial_q V(q) + \sum_{n=1}^{\ell} \delta^n \frac{1}{n!} \sum_{m=0}^{n\ell} \delta^m \sum_{\substack{k_1 + ... + k_n = m, \\ 0 \leq k_i \leq \ell}} \partial_q^{n+1} V(q)(g_{k_1,\ell}(q, p), ..., g_{k_n,\ell}(q, p)) + \delta^{\ell+1} \theta_\ell(q, p)$$

$$= \partial_q V(q) + I_{1,\ell}(q, p) + \delta^{\ell+1} I_{2,\ell}(q, p) + \delta^{\ell+1} \theta_\ell(q, p),$$

where

$$\theta_\ell(q, p) = \int_0^1 \frac{(1-t)^\ell}{\ell!} \partial_q^{\ell+2} V(q + t\delta z_1)(z_1, ..., z_1) dt,$$

$$I_{1,\ell}(q, p) = \sum_{n=1}^{\ell} \frac{1}{n!} \delta^n \sum_{m=0}^{j-n} \delta^m \sum_{\substack{k_1 + ... + k_n = m, \\ 0 \leq k_i \leq \ell}} \partial_q^{n+1} V(q)(g_{k_1,\ell}(q, p), ..., g_{k_n,\ell}(q, p)),$$

$$= \sum_{j=1}^{\ell} \delta^j \sum_{n=1}^{j} \frac{1}{n!} \sum_{\substack{k_1 + ... + k_n = \ell - n, \\ 0 \leq k_i \leq \ell - 1}} \partial_q^{n+1} V(q)(d_{k_1+1}(q, p), ..., d_{k_n+1}(q, p)),$$

$$I_{2,\ell}(q, p) = \sum_{n=1}^{\ell} \frac{1}{n!} \delta^n \sum_{m=\ell-n+1}^{n\ell} \delta^{m-\ell-1} \sum_{\substack{k_1 + ... + k_n = m, \\ 0 \leq k_i \leq \ell}} \partial_q^{n+1} V(q)(g_{k_1,\ell}(q, p), ..., g_{k_n,\ell}(q, p)).$$

Let $N$ be fixed, using the above computations, we have

$$q_1 = q + (1 + \gamma\delta)^{-1}(\delta p - \delta^2 \partial_q V(q_1))$$

$$= q + \sum_{k=0}^{N-1} (-\gamma)^k \delta^{k+1} p + \delta^{N+1} g(q, p)$$

$$+ \sum_{k=0}^{N-2} (-\gamma)^k \delta^{k+2} \Big( \partial_q V(q) - \delta^2 I_{1,N-k-2}(q, p) - \delta^{N-k-1} \big( I_{2,N-k-2}(q, p) + \theta_{N-k-2}(q, p) \big) \Big)$$

$$= q + \delta p + \sum_{k=2}^{N} (-1)^{k-1} \gamma^{k-2} \delta^k (\gamma p + \partial_q V(q)) + J(q, p) + \delta^{N+1} G(q, p), \tag{3.4}$$

where

$$J(q,p) = -\sum_{k=0}^{N-2}(-\gamma)^k I_{1,N-k-2}(q,p)\delta^{k+2}$$
$$= \sum_{k=2}^{N}(-1)^{k-1}\gamma^{k-2}\sum_{j=1}^{N-k}\delta^{j+k}\sum_{n=1}^{j}\frac{1}{n!}\sum_{\substack{k_1+\ldots+k_n=j-n,\\0\leq k_i\leq N}}\partial_q^{n+1}V(q)(d_{k_1+1}(q,p),\ldots,d_{k_n+1}(q,p))$$
$$= \sum_{k=3}^{N}\delta^k\sum_{j=2}^{k-1}(-1)^{j-1}\gamma^{j-2}\sum_{n=1}^{k-j}\frac{1}{n!}\sum_{\substack{k_1+\ldots+k_n=k-j-n,\\0\leq k_i\leq N}}\partial_q^{n+1}V(q)(d_{k_1+1}(q,p),\ldots,d_{k_n+1}(q,p)),$$
$$G(q,p) = g(q,p) - \sum_{k=0}^{N-1}(-\gamma)^k I_{2,N-k-2}(q,p),$$
$$g(q,p) = -\frac{(-\gamma)^{N-1}}{1+\delta\gamma}(\gamma p + \partial_q V(q_1)).$$

By identifying of (3.4) and (3.2), we get

$$q_1 = q + \sum_{k=1}^{N}\delta^k d_k(q,p) + \delta^{N+1}R_{N+1}(\delta,q,p),$$

where, for all $k \geq 0$, $d_k$ is defined for all $(q,p) \in \mathbb{R}^{2d}$ by

$$d_1(q,p) = p, \text{ and for } k \geq 2,$$
$$d_k(q,p) = (-1)^{k-1}\gamma^{k-2}(\gamma p + \partial_q V(q))$$
$$+ \sum_{j=2}^{k-1}(-1)^{j-1}\gamma^{j-2}\sum_{n=1}^{k-j}\frac{1}{n!}\sum_{\substack{k_1+\ldots+k_n=k-j-n,\\0\leq k_i}}\partial_q^{n+1}V(q)(d_{k_1+1}(q,p),\ldots,d_{k_n+1}(q,p)),$$

moreover, by induction, we have, for all $k \in \mathbb{N}$, $d_k \in \mathscr{C}_{pol}^{\infty}(\mathbb{R}^{2d})$.

The above identifying does not give an easy expression of $R_N$, then we have not immediately (3.7). To show this result, we will use that, for $N$, $q \in \mathbb{R}^d$ and $p \in \mathbb{R}^d$ fixed, $R_N(q,p,.)$ is the remainder of order $N$ of $\delta \mapsto \Psi_\delta(q,p)$. Therefore, if we show that, for any $n \in \mathbb{N}$, there exist $C > 0$ and $\ell_n \in \mathbb{N}$ such that for any $0 < \delta < \delta_0$ and $(q,p) \in \mathbb{R}^{2d}$,

$$|\partial_\delta^n \Psi_\delta(q,p)|^2 \leq C(1+|q|^{\ell_n}+|p|^{\ell_n}), \tag{3.5}$$

then we show (3.7) and Lemma 3.4 is shown.

Let us show the result (3.5) by induction on $n$. Let $0 < \delta < \delta_0$ and $(q,p) \in \mathbb{R}^{2d}$ be fixed. We have, by definition of $\Psi_\delta$,
$$\Psi_\delta(q,p)(1+\delta\gamma) + \delta^2\partial_q V(\Psi_\delta(q,p)) = \delta p + (1+\delta\gamma)q.$$

We multiply this equation by $\Psi_\delta(q,p)$ and use dissipativity inequality (2.3). We get

$$\langle p,\Psi_\delta(q,p)\rangle\delta + \langle q,\Psi_\delta(q,p)\rangle(1+\delta\gamma) = (1+\delta\gamma)|\Psi_\delta(q,p)|^2 + \delta^2\langle\partial_q V(\Psi_\delta(q,p)),\Psi_\delta(q,p)\rangle$$
$$\geq (1+\delta\gamma+\beta_1\delta^2)|\Psi_\delta(q,p)|^2 - \kappa\delta^2.$$

Since
$$2\langle a,b\rangle \leq \varepsilon|a|^2 + \frac{1}{\varepsilon}|b|^2 \quad \text{for all } \varepsilon > 0, a \in \mathbb{R}^d \text{ and } b \in \mathbb{R}^d, \tag{3.6}$$

we get, for any positive constants $\varepsilon_1$ and $\varepsilon_2$:

$$|\Psi_\delta(q,p)|^2 \left[1 + \delta\gamma + \beta_1\delta^2 - \frac{\varepsilon_1}{2} - \frac{\varepsilon_2}{2}\right]$$
$$\leq \kappa\delta^2 + \frac{\delta}{2\varepsilon_2}|p|^2 + \frac{1+\delta\gamma}{2\varepsilon_1}|q|^2.$$

We choose $\varepsilon_1 = \varepsilon_2 = \frac{1}{2}$, then we have that there exists a positive constant $C$, which is independent of $\delta$, such that

$$|\Psi_\delta(q,p)|^2 \leq C(1 + |q|^2 + |p|^2).$$

This proves (3.5) for $n = 0$ and $\ell_0 = 2$.

Let us assume the result (3.5) is true for all $0 \leq j < n$ and let us show it for $n$. Let $0 < \delta < \delta_0$ and $(q,p) \in \mathbb{R}^{2d}$ be fixed. We have

$$\begin{aligned}(1 + \delta\gamma)\partial_\delta^n \Psi_\delta(q,p) = &-n\gamma\partial_\delta^{n-1}\Psi_\delta(q,p) + (\gamma q + p)\mathbf{1}_{\{n=1\}}\\ &- 2n\partial_\delta^{n-1}\left(\partial_q V(\Psi_\delta(q,p))\right) - \delta^2 \partial_\delta^n\left(\partial_q V(\Psi_\delta(q,p))\right)\\ &- 2\binom{n}{2}\delta\partial_\delta^{n-2}\left(\partial_q V(\psi_\delta(q,p))\right)\\ =:&B_1(\delta,q,p) - 2nB_2(\delta,q,p) - \delta^2 B_3(\delta,q,p) - \delta B_4(\delta,q,p).\end{aligned}$$

By induction hypothesis, we have that $B_1$ has polynomial growth in $(q,p)$. Moreover, using Faà di Bruno's formula, we get

$$B_2(\delta,q,p) = \sum \frac{(n-1)!}{m_1! m_2!(2!)^{m_2}\dots m_{n-1}!((n-1)!)^{m_{n-1}}} \partial_q^{m_1+\dots+m_{n-1}+1} V(\Psi_\delta(q,p))$$
$$\times \prod_{j=1}^{n-1} \left(\partial_\delta^j \Psi_\delta(q,p)\right)^{m_j},$$

where $m_1 + 2m_2 + \dots + (n-1)m_{n-1} = n-1$. Using the polynomial growth of $V$ and its derivatives (**B-3**) and induction hypothesis, we have that $B_2$ has polynomial growth in $(q,p)$. Using the same method, we have that $B_4$ has polynomial growth in $(q,p)$. Using Faà di Bruno's formula, we also have

$$B_3(\delta,q,p) = \sum \frac{n!}{m_1! m_2!(2!)^{m_2}\dots m_n!(n!)^{m_n}} \partial_q^{m_1+\dots+m_n+1} V(\Psi_\delta(q,p))$$
$$\times \prod_{j=1}^{n} \left(\partial_\delta^j \Psi_\delta(q,p)\right)^{m_j}$$
$$= \delta \partial_q^2 V(q + \delta\Psi_\delta(q,p))\partial_\delta^n \Psi_\delta(q,p) + B_5(\delta,q,p),$$

where $m_1 + 2m_2 + \dots + nm_n = n$ and

$$B_5(\delta,q,p) = \sum \frac{n!}{k_1! k_2!(2!)^{k_2}\dots k_{n-1}!((n-1)!)^{k_{n-1}}} \partial_q^{k_1+\dots+k_{n-1}+1} V(\Psi_\delta(q,p))$$
$$\times \prod_{j=1}^{n-1} \left(\partial_\delta^j \Psi_\delta(q,p)\right)^{k_j},$$
$$k_1 + 2k_2 + \dots + (n-1)k_{n-1} = n.$$

Then we have that $B_1$, $B_2$, $B_4$ and $B_5$ have polynomial growth in $(q,p)$ and we get

$$\left((1+\delta\gamma)I + \delta^2 \partial_q^2 V(q + \delta\Psi_\delta(q,p))\right)\partial_\delta^n \Psi_\delta(q,p) = B_6(\delta,q,p)$$

where $B_6 := -n\gamma B_1 - 2nB_2 - \delta B_4 - \delta^2 B_5$ has polynomial growth in $(q,p)$.

Multiplying by $\delta^n_\delta \Psi_\delta(q,p)$, and using (3.6) on $|\langle B_6(\delta,q,p), \partial^n_\delta \Psi_\delta(q,p)\rangle|$ and semi-convexity assumption **B-1**, we get

$$|\partial^n_\delta \Psi_\delta(q,p)|^2 (1 + \gamma\delta - \theta\delta^2 - \frac{1}{2}) \leq \frac{1}{2} B_6^2(\delta,q,p),$$

where $\theta$ is the constant of semi-convexity and $B_6$ has polynomial growth in $(q,p)$. Since $\delta < \frac{\gamma}{\theta}$, we have that there exist constants $C$ and $\ell_n \in \mathbb{N}$ such that

$$|\partial^n_\delta \Psi_\delta(q,p)|^2 \leq C(1 + |q|^{\ell_n} + |p|^{\ell_n}),$$

which proves (3.5).

**Corollary 3.5** *Let $\delta_0 = \min(\frac{1}{\gamma}, \frac{\gamma\beta}{4\theta})$ and $(q^0, p^0) \in \mathbb{R}^{2d}$ such that $q_0 = q^0$ and $p_0 = p^0$. We have for $0 < \delta < \delta_0$ and $N \in \mathbb{N}$,*

$$q_1 = q^0 + \sum_{k=1}^{N} \delta^k d_k(q^0, p^0) + \delta^{N+1} R_{N+1}(q^0, p^0, \delta),$$

$$p_1 = \sum_{k=0}^{N} \delta^k d_{k+1}(q^0, p^0) + \delta^{N+1} R_{N+2}(q^0, p^0, \delta) + \sqrt{\delta}\sigma\eta_0$$

*where, for all $k \geq 0$, $d_k \in \mathscr{C}^\infty_{pol}(\mathbb{R}^{2d})$ is defined for all $(x,y) \in \mathbb{R}^{2d}$ by*

$$d_1(x,y) = y, \text{ and for } k \geq 2$$
$$d_k(x,y) = (-1)^{k-1}\gamma^{k-2}(\gamma y + \partial_q V(x))$$
$$+ \sum_{j=2}^{k-1}(-1)^{j-1}\gamma^{j-2}\sum_{n=1}^{k-j}\frac{1}{n!}\sum_{\substack{k_1+\ldots+k_n=k-j-n,\\ 0 \leq k_i}} \partial^{n+1}_q V(x)(d_{k_1+1}(x,y), \ldots, d_{k_n+1}(x,y)).$$

*and $R_{N+1}$ verifies: There exist $C > 0$ and $\ell_N \in \mathbb{N}$ such that for any $(x,y) \in \mathbb{R}^{2d}$ and $\delta < \delta_0$*

$$|R_{N+1}(x,y,\delta)| \leq C(1 + |x|^{\ell_N} + |y|^{\ell_N}). \tag{3.7}$$

*Proof (Proof of Proposition 3.2 in the case of the implicit split-step scheme (2.11).)* Let $N$ fixed. We have, with the notation of Corollary 3.5, for all $N_0 \in \mathbb{N}$

$$q_1 = q + \sum_{k=1}^{N_0} \delta^k d_k(q,p) + \delta^{N_0+1} R_{N_0+1}(q,p,\delta) = q + \delta R_1(q,p,\delta),$$

$$p_1 = \sum_{k=0}^{N_0} \delta^k d_{k+1}(q,p) + \delta^{N_0+1} R_{N_0+2}(q,p,\delta) + \sqrt{\delta}\sigma\eta_0 = p + \sqrt{\delta}\sigma\eta_0 + \delta R_2(q,p,\delta) = z + \sqrt{\delta}\sigma\eta_0.$$

Let $\phi \in \mathscr{C}^\infty_{pol}(\mathbb{R}^{2d})$, for any $n \in \mathbb{N}$, there exists an integer $r_n$ such that $\phi \in \mathscr{C}^n_{r_n}(\mathbb{R}^{2d})$. Using Taylor expansion, we get

$$\phi(q_1, p_1) = \phi(q_1, z + \sqrt{\delta}\sigma\eta_0)$$
$$= \phi(q_1, z) + \sum_{k=1}^{2N+1} \frac{1}{k!} \delta^{k/2} \sigma^k \partial^k_p \phi(q_1, z)(\eta_0, \ldots \eta_0)$$
$$+ \int_0^1 \frac{(1-t)^{2N+1}}{(2N+1)!} \delta^{N+1}\sigma^{2N+2} \partial^{2N+2}_p \phi(q_1, z + t\sqrt{\delta}\sigma\eta_0)(\eta_0, \ldots \eta_0) dt.$$

Let $\lfloor . \rfloor$ denotes the integer part. Using Taylor expansion on $\partial_p^k \phi(q_1, z)$ and computations done in the proof of Lemma 3.4, we obtain

$$\phi(q_1, p_1) = I_1(x, \eta_0) + I_2(x, \eta_0) + I_3(x, \eta_0) + I_4(x, \eta_0),$$

where $x = (q, p)$ and

$$I_1(x, \eta_0) = \sum_{k=0}^{2N+1} \frac{1}{k!} \sigma^k \Big( \delta^{k/2} \partial_p^k \phi(x)(\eta_0, ..., \eta_0) + \sum_{n=1}^{\ell_k} \frac{1}{n!} \sum_{\ell=0}^{n} \binom{n}{\ell} \sum_{m=n}^{\ell_k} \delta^{m+k/2}$$
$$\times \sum_{\substack{k_1+...+k_\ell \\ +\tilde{k}_1+...+\tilde{k}_{n-\ell}=m \\ 0<k_i \leq \ell_k \\ 0<\tilde{k}_i \leq \ell_k}} \partial_q^\ell \partial_p^{n-\ell+k} \phi(x)(d_{k_1}(x), ..., d_{k_\ell}(x), d_{\tilde{k}_1+1}(x), ..., d_{\tilde{k}_{n-\ell}+1}(x), \eta_0, ..., \eta_0) \Big),$$

$$I_2(x, \eta_0) = \sum_{k=0}^{2N+1} \frac{1}{k!} \sigma^k \sum_{n=0}^{\ell_k} \frac{1}{n!} \sum_{\ell=0}^{n} \binom{n}{\ell} \sum_{m=\ell_k+1}^{n(\ell_k+1)} \delta^{m+k/2} B_{m,n,\ell,k}(x_1, \eta_0),$$

$$I_3(x, \eta_0) = \sum_{k=0}^{2N+1} \frac{1}{k!} \sigma^k \delta^{k/2} \delta^{\ell_k+1} \int_0^1 \frac{(1-t)^{\ell_k}}{(l_k)!} \sum_{\ell=0}^{\ell_k+1} \binom{\ell_k+1}{\ell} \partial_q^\ell \partial_p^{\ell_k+1-\ell+k} \phi(q+t\delta R_1(x,\delta), p+t\delta R_2(x,\delta))$$
$$\times (R_1(x,\delta), ..., R_1(x,\delta), R_2(x,\delta), ..., R_2(x,\delta), \eta_0, ..., \eta_0) dt,$$

$$I_4(x, \eta_0) = \delta^{N+1} \sigma^{2N+2} \int_0^1 \frac{(1-t)^{2N+1}}{(2N+1)!} \partial_p^{2N+2} \phi(q+\delta R_1(x,\delta), z+t\sqrt{\delta}\eta_0)(\eta_0, ..., \eta_0) dt,$$

$$\ell_k = N - \left\lfloor \frac{k+1}{2} \right\rfloor$$

and, with the temporary notation: for all $0 \leq k \leq 2N+1$ and $1 \leq i \leq \ell_k$, $\tilde{g}_{k,i} = d_{i+1}$, $g_{k,i} = d_i$, $\tilde{g}_{k,\ell_k+1} = R_{\ell_k+2}$ and $g_{k,\ell_k+1} = R_{\ell_k+1}$,

$$B_{m,n,\ell,k}(x, \eta_0) = \sum_{\substack{k_1+...+k_\ell \\ +\tilde{k}_1+...+\tilde{k}_{n-\ell}=m \\ 0<k_i \leq \ell_k+1 \\ 0<\tilde{k}_i \leq \ell_k+1}} \partial_q^\ell \partial_p^{n-\ell+k} \phi(x)(g_{k_1}(x), ..., g_{k_\ell}(x), \tilde{g}_{\tilde{k}_1}(x), ..., \tilde{g}_{\tilde{k}_{n-\ell}}(x), \eta_0, ..., \eta_0).$$

We have, for all $\ell \in \mathbb{N}$, $\mathbb{E}(\eta_{0,i}^{2\ell+1}) = 0$ and $\eta_{0,i}$ are independent of $q$, $p$ and $\eta_{0,j}$ for $j \neq i$, then the expectation of all the odd term in $k$ in $I_1$, $I_2$ and $I_3$ vanish. Hence we have

$$\mathbb{E}(I_2(x, \eta_0)) = \sum_{k=0}^{N} \frac{1}{(2k)!} \sigma^{2k} \sum_{n=0}^{N-k} \frac{1}{n!} \sum_{\ell=0}^{n} \binom{n}{\ell} \sum_{m=N-k+1}^{n(N-k+1)} \delta^{m+k} \mathbb{E}(B_{m,n,\ell,2k}(x, \eta_0)),$$

$$\mathbb{E}(I_3(x, \eta_0)) = \sum_{k=0}^{N} \frac{1}{(2k)!} \delta^{N+1} \sigma^{2k} \int_0^1 \frac{(1-t)^{N-k}}{(N-k)!} \sum_{\ell=0}^{N-k+1} \binom{N-k+1}{\ell}$$
$$\times \mathbb{E}\Big( \partial_q^\ell \partial_p^{N-k+1-\ell+k} \phi(q+t\delta R_1(x,\delta), p+t\delta R_2(x,\delta))$$
$$\times (R_1(x,\delta), ..., R_1(x,\delta), R_2(x,\delta), ..., R_2(x,\delta), \eta_0, ..., \eta_0) \Big) dt.$$

In $\mathbb{E}(I_2(x, \eta_0))$, we have $m+k \geq N+1$, then we can factor $\delta^{N+1}$. Moreover, if $n = \ell = 0$ then $B_{m,n,\ell,2k} = 0$, hence, in each term of $\mathbb{E}(I_2(x, \eta_0))$, we have at least one derivative of $\phi$. Using $\phi \in$ $\mathscr{C}^{2N+2}_{r_{2N+2}}(\mathbb{R}^{2d})$ and the polynomial growth of $d_j$, $\tilde{R}_{N_0+1}$ and $R_{N_0+1}$ for $j \in \mathbb{N}^*$ and $N_0 \in \mathbb{N}$, we get that there exist integers $n_1$, $n_2$ and $n_3$ such that

$$|\mathbb{E}(I_2(q,p,\eta_0))| \leq C_N \delta^{N+1}(1+|p|^{n_1}+|q|^{n_1})|\phi|_{2N,r_{2N}},$$
$$|\mathbb{E}(I_4(q,p,\eta_0))| \leq C_N \delta^{N+1}(1+|p|^{n_2}+|q|^{n_2}) \|D^{2N+2}\phi\|_{0,r_{2N+2}}$$
$$|\mathbb{E}(I_3(q,p,\eta_0))| \leq C_N \delta^{N+1}(1+|p|^{n_3}+|q|^{n_3}) \|D^N \phi\|_{N+1,r_{2N+1}}.$$

Hence, we have that there exists $k_1 \in \mathbb{N}^*$ such that

$$|\mathbb{E}\phi(q_1,p_1) - \mathbb{E}(I_1)| \leq C_N \delta^{N+1}(1+|p|^{k_1}+|q|^{k_1})|\phi|_{2N+2,r_{2N+2}},$$

where

$$\mathbb{E}(I_1) = \sum_{j=0}^{N} \frac{1}{(2j)!} \sigma^{2j} \delta^j \mathbb{E}\left(\partial_p^{2j}\phi(q,p)(\eta_0,...,\eta_0)\right)$$
$$+ \sum_{j=0}^{N} \frac{1}{(2j)!} \sigma^{2j} \sum_{n=1}^{N-j} \frac{1}{n!} \sum_{\ell=0}^{n} \binom{n}{\ell} \sum_{m=n}^{N-j} \delta^{m+j} \sum_{\substack{k_1+...+k_\ell \\ +\tilde{k}_1+...+\tilde{k}_{n-\ell}=m \\ \tilde{k}_i, k_i \geq 1}} \mathbb{E}\left(\partial_q^\ell \partial_p^{n-\ell+2j}\phi(q,p)\right.$$
$$\left. \times (d_{k_1}(q,p),...,d_{k_\ell}(q,p),d_{\tilde{k}_1+1}(q,p),...,d_{\tilde{k}_{n-\ell}+1}(q,p),\eta_0,...,\eta_0)\right),$$
$$= \sum_{k=0}^{N} \delta^k A_k(q,p)\phi(q,p),$$

and

$$A_k(q,p)\phi(q,p) = \frac{\sigma^{2k}}{(2k)!}\mathbb{E}\left(\partial_p^{2k}\phi(q,p)(\eta_0,...,\eta_0)\right) + \sum_{j=0}^{k-1} \frac{1}{(2j)!}\sigma^{2j} \sum_{n=1}^{k-j}\frac{1}{n!}\sum_{\ell=0}^{n}\binom{n}{\ell}$$
$$\times \sum_{\substack{k_1+...+k_\ell \\ +\tilde{k}_1+...+\tilde{k}_{n-\ell}=k-j \\ 0<k_i,\tilde{k}_i}} \mathbb{E}\left(\partial_q^l \partial_p^{n-\ell+2j}\phi(q,p)(d_{k_1}(q,p),...,d_{k_\ell}(q,p),d_{\tilde{k}_1+1}(q,p),...,d_{\tilde{k}_{n-\ell}+1}(q,p),\eta_0,...,\eta_0)\right)$$

Using property of $\eta_0$, we get that $A_0 = I$,

$$A_1\phi(q,p) = \frac{\sigma^2}{2}\sum_{i=1}^{d} \frac{\partial^2}{\partial_{p_i}\partial_{p_i}}\phi(q,p) + \langle \partial_p\phi, d_2(q,p) \rangle + \langle \partial_q\phi(q,p), d_1(q,p) \rangle = L$$

and $A_k$ is an operator of order $2k$.

The proof in the case of the scheme (2.12) uses the same arguments.
We need asymptotic expansions for $q_1 = q + \delta p_1$ and $p_1 = p - \delta \gamma p_1 - \partial_q V(q_1)\delta + \sqrt{\delta}\sigma\eta_0$. We use the local notation $\alpha = \sqrt{\delta}$. We define the function $\psi_\alpha$ which associate to $(q,p)$ the solution $z$ of $(1+\gamma\alpha^2)z = (1+\gamma\alpha^2)q + \alpha^2 p - \alpha^4 \partial_q V(z) + \alpha\sigma\eta_0$. This function is well defined (see Lemma 2.7). Moreover, using same arguments as for the scheme (2.11), we can show that $(\alpha,q,p) \mapsto \psi_\alpha(q,p)$ is $C^\infty$ on $]0,\sqrt{\delta_1}[\times \mathbb{R}^{2d}$, where $\delta_1 = \frac{\gamma+\sqrt{\gamma^2+4\theta}}{2\theta}$.

We have the following lemma:

**Lemma 3.6** *Let $\delta_0 = \min(\frac{1}{\gamma}, \frac{\gamma\beta}{4\theta})$. For $0 < \delta < \delta_0$ and the local notation $\alpha = \sqrt{\delta}$, we have*

$$\forall N_0 \in \mathbb{N}, \quad q_1 = \psi_\alpha(q,p) = q + \sum_{k=2}^{2N_0+1} \delta^{\frac{k}{2}} d_k(q,p,\eta_0) + \delta^{N_0+1} R_{N_0+1}(q,p,\delta,\eta_0),$$

$$p_1 = \sum_{k=0}^{2N_0+1} \delta^{\frac{k}{2}} d_{k+2}(q,p,\eta_0) + \delta^{N_0+1} R_{N_0+2}(q,p,\delta,\eta_0),$$

*where, $\forall k \geq 2$, $d_k$ is defined for all $(q,p) \in \mathbb{R}^{2d}$ by*

$$d_2(q,p,\eta_0) = p, d_3(q,p,\eta_0) = \sigma\eta_0$$

*and $\forall k \geq 2$*

$$d_{2k}(q,p,\eta_0) = (-1)^{k-1}\gamma^{k-2}(\gamma p + \partial_q V(q))$$
$$+ \sum_{i=2}^{k-1} (-1)^{i+1}\gamma^i \sum_{n=1}^{k-i} \frac{1}{n!} \sum_{\substack{k_1+\ldots+k_n=2(k-i) \\ k_s \geq 2}} \partial_q^n V(q)(d_{k_1}(q,p,\eta_0),\ldots,d_{k_i}(q,p,\eta_0)),$$

$$d_{2k+1}(q,p,\eta_0) = (-1)^{k-1}\gamma^{k-1}\sigma\eta_0$$
$$+ \sum_{i=2}^{k} (-1)^{i+1}\gamma^i \sum_{n=1}^{k-i} \frac{1}{n!} \sum_{\substack{k_1+\ldots+k_n=2(k-i)+1 \\ k_s \geq 2}} \partial_q^n V(q)(d_{k_1}(q,p,\eta_0),\ldots,d_{k_i}(q,p,\eta_0)).$$

*Moreover, we have that for any $k \geq 2$, $\mathbb{E}(d_k) \in \mathscr{C}_{pol}^\infty(\mathbb{R}^{2d})$,*

$$\mathbb{E}(d_{2k+1}(q,p,\eta_0)) = 0 \tag{3.8}$$

*and, for any $N \in \mathbb{N}$, $R_N$ verifies: There exist $C > 0$ and $\ell_N \in \mathbb{N}$ such that for any $(x,y) \in \mathbb{R}^{2d}$ and $\delta < \delta_0$*

$$|\mathbb{E}(R_N(q,p,\delta,\eta_0))| \leq C(1 + |q|^{\ell_N} + |p|^{\ell_N}).$$

*Proof* To prove this Lemma, we use the same ideas as in Lemma 3.4. We first compute $d_k$ for all $k$. By induction, we rewrite $d_k$ only in terms of $d_1$, $d_2$ and the derivatives of $V$ evaluate in $q$. Using the independence of $\eta_0$ with $(q,p)$, we can show (3.8).
To prove that $\mathbb{E}(R_N)$ has polynomial growth, we show that, for $n \in \mathbb{N}$, there exist $C_n > 0$ and $k_n \in \mathbb{N}$ such that, for $(q,p) \in \mathbb{R}^{2d}$, $\delta < \delta_0$ and the local notation $\alpha = \sqrt{\delta}$, we have

$$\mathbb{E}(|\partial_\alpha^n \psi_\alpha(q,p)|^2) \leq C_n(1 + |q|^{k_n} + |p|^{k_n}).$$

The proof of Proposition 3.2 in the case of the scheme (2.12) is similar to the case of the scheme (2.11), but we must use an asymptotic expansion of $\partial_p^k \phi$ to a larger order ($2N+1-k$ instead of $N - \lfloor (k+1)/2 \rfloor$).

## 4 Modified generator

For now on, all definitions depend of the scheme considered.

## 4.1 Formal series analysis

Let us now consider $\delta$ as fixed. We want to construct a formal series

$$\mathscr{L} = L + \delta L_1 + ... + \delta^n L_n + ... \tag{4.1}$$

where the coefficients of the operator $L_n$ are in $\mathscr{C}^\infty_{pol}(\mathbb{R}^{2d})$ and such that formally the solution $v$ at time $t = \delta$ of the equation

$$\partial_t v(t,q) = \mathscr{L}(q,p)v(t,q,p), t > 0, (q,p) \in \mathbb{R}^{2d} \quad v(0,q,p) = \phi(q,p), (q,p) \in \mathbb{R}^{2d}$$

coincides in the sense of asymptotic expansion with the approximation of the transition semigroup $\mathbb{E}\phi(q_1, p_1)$ studied in the previous section. In other words, we want to have the equality in the sense of asymptotic expansion in powers of $\delta$

$$\exp(\delta\mathscr{L})\phi = \phi + \sum_{n \geq 1} \delta^n A_n \phi,$$

where the operators $A_n$ are defined in Proposition 3.2.

Formally, this equation can be written as

$$\exp(\delta\mathscr{L}) - I_d = \delta \tilde{A}(\delta), \tag{4.2}$$

where $\tilde{A}(\delta) = \sum_{n \geq 1} \delta^{n-1} A_n$.

We have

$$\exp(\delta\mathscr{L}) - I_d = \delta\mathscr{L}\left(\sum_{n \geq 0} \frac{\delta^n}{(n+1)!}\mathscr{L}^n\right).$$

Note that the (formal) inverse of the series is given by

$$\left(\sum_{n \geq 0} \frac{\delta^n}{(n+1)!}\mathscr{L}^n\right)^{-1} = \sum_{n \geq 0} \frac{B_n}{n!}\delta^n \mathscr{L}^n,$$

where the $B_n$ are the Bernoulli numbers (see [6,9]). Hence, equations (4.1) and (4.2) are equivalent in the sense of formal series to

$$\mathscr{L} = \sum_{\ell \geq 0} \frac{B_\ell}{\ell!}\delta^\ell \mathscr{L}^\ell \tilde{A}(\delta) = \sum_{n \geq 0} \delta^n \left(A_{n+1} + \sum_{\ell=1}^n \frac{B_\ell}{\ell!} \sum_{n_1+...+n_{\ell+1}=n-\ell} L_{n_1}...L_{n_\ell} A_{n_{\ell+1}+1}\right). \tag{4.3}$$

Identifying the right hand sides of (4.1) and (4.3), we get the following induction formula

$$L_n = A_{n+1} + \sum_{\ell=1}^n \frac{B_\ell}{\ell!} \sum_{n_1+...+n_{\ell+1}=n-\ell} L_{n_1}...L_{n_\ell} A_{n_{\ell+1}+1}. \tag{4.4}$$

Each of the terms of the above sum is an operator of order $2n+2$ with coefficients $\mathscr{C}^\infty_{pol}(\mathbb{R}^{2d})$ and therefore $L_n$ is also an operator of order $2n+2$ with coefficients $\mathscr{C}^\infty_{pol}(\mathbb{R}^{2d})$.

Notes that (4.2) gives immediately the inverse relation of this formal series equation:

$$A_n = \sum_{\ell=1}^n \frac{1}{\ell!} \sum_{n_1+...+n_\ell=n-\ell} L_{n_1}...L_{n_\ell}. \tag{4.5}$$

Moreover we have clearly

$$L_n 1 = 0.$$

4.2 Approximate solution of the modified flow

For a given $N$, we have constructed in the previous section a modified operator

$$L^{(N)} = L + \sum_{n=1}^{N} \delta^n L_n. \tag{4.6}$$

In order to perform weak backward error analysis and estimate recursively the modified invariant law of the numerical process, we should be able to define a solution $v^N$ of the modified flow

$$\partial_t v^N(t,q,p) = L^{(N)} v^N(t,q,p), t > 0, (q,p) \in \mathbb{R}^{2d} \quad v^N(0,q,p) = \phi(q,p), (q,p) \in \mathbb{R}^{2d}. \tag{4.7}$$

However in our situation we do not know whether this equation has a solution.

The goal of the following theorem is to give a proper definition of a modified flow associated to (4.6).

**Theorem 4.1** *Let $\phi \in \mathscr{C}_{pol}^{\infty}(\mathbb{R}^{2d})$ such that $\int \phi d\rho = 0$. For any $N \in \mathbb{N}$, there exists an integer $r_N$ such that $\phi \in \mathscr{C}_{r_N}^{N}(\mathbb{R}^{2d})$. For all $n \in \mathbb{N}$, there exist functions $v_n(t,.) \in \mathscr{C}_{pol}^{\infty}(\mathbb{R}^{2d})$ defined for all times $t \geq 0$ such that for all $t \geq 0$ and $n \in \mathbb{N}$,*

$$\partial_t v_n(t,q,p) - L v_n(t,q,p) = \sum_{\ell=1}^{n} L_\ell v_{n-\ell}(t,q,p), \tag{4.8}$$

*with initial condition $v_0(0,.) = \phi$ and $v_n(0,.) = 0$ for $n \geq 1$. For all $N \geq 0$, setting*

$$v^{(N)}(t,q,p) = \sum_{k=0}^{N} \delta^k v_k(t,q,p),$$

*then the following holds:*

a. *Let $\delta_0 = \min(\frac{1}{\gamma}, \frac{\gamma\beta}{4\theta})$, there exist a positive real number $C_N$ and integers $\ell_N$ and $k_N$ such that for all $t \geq 0$, $0 < \delta < \delta_0$ and $(q,p) \in \mathbb{R}^{2d}$*

$$|\mathbb{E} v^{(N)}(t,q_1,p_1) - v^{(N)}(t+\delta,q,p)| \leq \delta^{N+1} C_N (1+|q|^{\ell_N}+|p|^{\ell_N}) \sup_{\substack{s \in ]0,\delta[ \\ n=0,\ldots,N}} |v_n(t+s,.)|_{2N+2,k_N}.$$

b. *Let $\delta_0 = \min(\frac{1}{\gamma}, \frac{\gamma\beta}{4\theta})$, there exist a positive real number $C_N$ and an integer $\ell_N$ such that for all $0 < \delta < \delta_0$ and $(q,p) \in \mathbb{R}^{2d}$*

$$|\mathbb{E}\phi(q_1,p_1) - v^{(N)}(\delta,q,p)| \leq \delta^{N+1} C_N (1+|q|^{\ell_N}+|p|^{\ell_N}) \| \phi \|_{N(d+7)+2, r_{N(d+7)+2}}.$$

*Proof* Let $\phi \in \mathscr{C}_{pol}^{\infty}(\mathbb{R}^{2d})$ such that $\int \phi d\rho = 0$, for any $N \in \mathbb{N}$, there exists an integer $r_N$ such that $\phi \in \mathscr{C}_{r_N}^{N}(\mathbb{R}^{2d})$. For $n = 0$, equation (4.8) reduces to $v_0 = u$, the solution of (2.27). By Proposition 2.15, we have that $u$ and all its derivatives have polynomial growth. Let $n \in \mathbb{N}$ and assume that $v_j(t)$ are constructed for $j = 1, ..., n-1$. Let for $t \geq 0$ and $(q,p) \in \mathbb{R}^{2d}$

$$F_n(t,q,p) = \sum_{l=1}^{n} L_l v_{n-l}(t,q,p), \tag{4.9}$$

the right-hand side in (4.8), then $v_n$ is uniquely defined and given by the formula

$$v_n(t,.) = \int_0^t P_{t-s} F_n(s,.) ds, \quad t \geq 0. \tag{4.10}$$

Using an induction argument and Proposition 2.15, we know that $v_n$ and all its derivatives have polynomial growth. Moreover, we have for all $k \in \mathbb{N}$ and $n \in \mathbb{N}^*$ that there exists an integer $\alpha_{k,n}$ such that for all $t \geq 0$

$$\| v_n(t) \|_{k,\alpha_{k,n}} \leq C(t) \| \phi \|_{k+n(d+1)+4n, r_{k+n(d+5)}}, \tag{4.11}$$

where the constant $C(t)$ depends on $t$, $k$, $n$ and $V$. This proves the first part of the Theorem.

To prove **a.**, we consider a fixed time t and define the functions $w_n(s,q,p) := v_n(t+s,q,p)$ for $s \geq 0$, $(q,p) \in \mathbb{R}^{2d}$ and $n \in \mathbb{N}$. By definition, these functions satisfy the relation

$$\partial_s w_n(s,q,p) = \sum_{\ell=0}^n L_\ell w_{n-\ell}(s,q,p), s \geq 0, (q,p) \in \mathbb{R}^{2d} \quad w_n(0,q,p) = v_n(t,q,p), (q,p) \in \mathbb{R}^{2d}.$$

Let us consider the successive time derivatives of the functions $w_n$. We have, using the definition of $w_n$, for all $s \geq 0$ and $(q,p) \in \mathbb{R}^{2d}$,

$$\partial_s^2 w_n(s,q,p) = \sum_{\ell=0}^n L_\ell \partial_s w_{n-\ell}(s,q,p) = \sum_{k=0}^n \sum_{\ell_1+\ell_2=k} L_{\ell_1} L_{\ell_2} w_{n-k}(s,q,p),$$

and we see by induction that for all $m \geq 1$, $(q,p) \in \mathbb{R}^{2d}$ and $s \geq 0$

$$\partial_s^m w_n(s,q,p) = \sum_{\ell_1+\ldots+\ell_{m+1}=n} L_{\ell_1}\ldots L_{\ell_m} w_{\ell_{m+1}}(s,q,p).$$

Using the fact that the operators $L_\ell$ are of order $2\ell+2$ with no terms of order zero and that their coefficients have polynomial growth, we see that there exist a constant $C$, depending on $n$ and $m$, and an integer $\beta_n$ such that

$$|\partial_s^m w_n(s,q,p)| \leq C(1+|q|^{\beta_n}+|p|^{\beta_n}) \sup_{k=0,\ldots,n} |w_k(s,.)|_{2n-2k+2m,\ell_{n,k,m}},$$

where $\ell_{n,k,m} \in \mathbb{N}$ is such that $w_k(s,.) \in \mathscr{C}^{2(n-k+m)}_{\ell_{n,k,m}}(\mathbb{R}^{2d})$. Now let us consider the Taylor expansion of $w_n(\delta,.)$, for $\delta < \delta_0$. Let $N$ be fixed. We have for $\delta < \delta_0$, $(q,p) \in \mathbb{R}^{2d}$ and $n = 0,\ldots,N$,

$$w_n(\delta,q,p) = \sum_{m=0}^{N-n} \frac{\delta^m}{m!} \partial_t^m w_n(0,q,p) + \int_0^\delta \frac{s^{N-n}}{(N-n)!} \partial_t^{N-n+1} w_n(s,q,p) ds$$

$$= \sum_{m=0}^{N-n} \frac{\delta^m}{m!} \sum_{\ell_1+\ldots+\ell_{m+1}=n} L_{\ell_1}\ldots L_{\ell_m} w_{\ell_{m+1}}(0,q,p) + R_{N,n}(\delta,q,p).$$

Using the bounds on the time derivatives of $w_n$, we obtain that there exists an integer $\ell_N$ such that for all $0 \leq \delta \leq \delta_0$ and all $n = 0,\ldots,N$,

$$|R_{N,n}(\delta,q,p)| \leq C\delta^{N-n+1}(1+|q|^{\ell_N}+|p|^{\ell_N}) \sup_{\substack{s\in ]0,\delta[,\\n=0,\ldots,N}} |w_n(s,.)|_{2N+2,k_N},$$

for some constants depending on $N$, $n$ and $k_N$ such that for $n = 0,\ldots,N$, $w_n \in \mathscr{C}^{2N+2}_{k_N}(\mathbb{R}^{2d})$.

After summation in $n$ and using the expression (4.5) of the operators $A_n$ and the definition of $w_n$, we get for all $(q,p) \in \mathbb{R}^{2d}$, $t \geq 0$ and $0 < \delta < \delta_0$

$$v^{(N)}(t+\delta,q,p) = \sum_{n=0}^N \delta^n \sum_{m=0}^n A_m v_{n-m}(t,q,p) + R_N(t,\delta,q,p),$$

where

$$|R_N(t,\delta,q,p)| \leq C_N \delta^{N+1}(1+|q|^{\ell_N}+|p|^{\ell_N}) \sup_{\substack{s\in ]0,\delta[,\\n=0,\ldots,N}} |v_n(t+s,.)|_{2N+2,k_N}.$$

To conclude we use (3.1) applied to $\phi = v^{(N)}(t,q,p)$ and the fact that $\delta < \delta_0$.
The second estimate **b.** is then a consequence of **a.** with $t = 0$ and (4.11).

## 5 Asymptotic expansion of the invariant measure and long time behavior

We now analyze the long time behavior of the solution of the modified flow associated to (4.7). In the following, for a given operator $B$, we denote by $B^*$ its adjoint with respect to the $L^2(\rho)$ product. We start by an asymptotic expansion of a formal invariant measure for the numerical schemes.

**Proposition 5.1** *Let $\delta_0 = \min(\frac{1}{\gamma}, \frac{\gamma\beta}{4\theta})$. Let $(L_n)_{n\geq 0}$ be the collection of operators defined recursively by (4.4). There exists a collection of functions $(\mu_n)_{n\geq 0}$ such that $\mu_0 = 1$ and $\int \mu_n d\rho = 0$ for all $n \geq 1$. Moreover, for all n, $\mu_n \in \mathscr{C}_{pol}^\infty(\mathbb{R}^{2d})$ and for all $n \geq 1$*

$$L^*\mu_n = -\sum_{\ell=1}^{n}(L_\ell)^*\mu_{n-\ell}. \tag{5.1}$$

*Let the function $\mu^{(N)}$ be defined for $(q,p) \in \mathbb{R}^{2d}$ and $0 < \delta < \delta_0$ by*

$$\mu^{(N)}(\delta,q,p) = 1 + \sum_{n=1}^{N} \delta^n \mu_n(q,p),$$

*then for $0 < \delta < \delta_0$, $\mu^{(N)}(\delta,.,.) \in \mathscr{C}_{pol}^\infty(\mathbb{R}^{2d})$ and satisfies*

$$\int_{\mathbb{R}^{2d}} \mu^{(N)}(\delta,q,p)\rho(q,p)dqdp = 1.$$

*Proof* Let $n \geq 1$ be fixed. Assume that $\mu_0 = 1$ and for $j = 1,...,n-1$, $\mu_j$ are known and $\mu_j \in \mathscr{C}_{pol}^\infty(\mathbb{R}^{2d})$. Let us consider equation (5.1) given by

$$L^*\mu_n = -\sum_{\ell=1}^{n} L_\ell^* \mu_{n-\ell} =: G_n$$

Note that $G_n \in \mathscr{C}_{pol}^\infty(\mathbb{R}^{2d})$. Indeed, $\mu_0, ..., \mu_{n-1}$ and all the coefficients of $L_\ell^*$ are in $\mathscr{C}_{pol}^\infty(\mathbb{R}^{2d})$. Moreover, $G_n$ satisfies

$$\int_{\mathbb{R}^{2d}} G_n(q,p)\rho(q,p)dqdp = -\sum_{\ell=1}^{n} \int_{\mathbb{R}^{2d}} L_\ell^*\mu_{n-\ell}(q,p)\rho(q,p)dqdp$$

$$= -\sum_{\ell=1}^{n} \int_{\mathbb{R}^{2d}} \mu_{n-\ell}(q,p)L_\ell 1 \rho(q,p)dqdp = 0.$$

Using Lemma 2.16, we easily obtain the existence of a function $\mu_n \in \mathscr{C}_{pol}^\infty(\mathbb{R}^{2d})$ satisfying (5.1) and $\int \mu_n d\rho = 0$. This shows the proposition. $\square$

**Proposition 5.2** *Let $\phi \in \mathscr{C}_{pol}^\infty(\mathbb{R}^{2d})$ such that $u(0) = \phi$, then for any $n \in \mathbb{N}$, there exists an integer $r_n$ such that $\phi \in \mathscr{C}_{r_n}^n(\mathbb{R}^{2d})$. For all $n \in \mathbb{N}$, $k \in \mathbb{N}$ and $\lambda < \lambda_0$, there exist a positive polynomial function $P_{k,n}$ and an integer $\ell_{k,n}$ such that for all $t \geq 0$*

$$\| v_n(t) - \int \phi \mu_n d\rho \|_{k,\ell_{k,n}} \leq P_{k,n}(t)e^{-\lambda t} \| \phi - \langle \phi \rangle \|_{\beta_{k,n},\alpha_{k,n}}, \tag{5.2}$$

*where $\lambda_0$ is defined in the Proposition 2.15, $\langle \phi \rangle = \int \phi \rho dqdp$, $\beta_{k,n} = k + 4n + (n+1)(d+1)$ and $\alpha_{k,n} = r_{\beta_{k,n}}$.*

*Proof* Let $\phi \in \mathscr{C}^{\infty}_{pol}(\mathbb{R}^{2d})$ such that $u(0) = \phi$, then for any $n \in \mathbb{N}$, there exists an integer $r_n$ such that $\phi \in \mathscr{C}^n_{r_n}(\mathbb{R}^{2d})$. Using the fact that $\mu_0 = 1$ and $v_0 = u$, we see that the estimate (5.2) is satisfied for $n = 0$ (Proposition 2.15). Let $n \geq 1$ assume that $v_j$, $j = 0,...,n-1$ satisfy: For all $k \in \mathbb{N}$ and $\lambda < \lambda_0$, there exist a positive polynomial function $P_{k,j}$ and $\ell_{k,l}$ such that for any $t \geq 0$:

$$\| v_j(t,.) - \int \phi \mu_j d\rho \|_{k,\ell_{j,k}} \leq P_{k,j}(t) e^{-\lambda t} \| \phi - \langle \phi \rangle \|_{k+4j+(j+1)(d+1), \alpha_{k,j}},$$

where $\alpha_{k,j}$ is such that $\phi \in \mathscr{C}^{k+4j+(j+1)(d+1)}_{\alpha_{k,j}}(\mathbb{R}^{2d})$. Let us set for $t \geq 0$

$$c_n(t) = \sum_{m=0}^{n} \int_{\mathbb{R}^{2d}} v_{n-m}(t,q,p)\mu_m(q,p)\rho(q,p)dqdp.$$

We claim that $c_n(.)$ does not depend on time. Indeed, for all $t \geq 0$,

$$\sum_{m=0}^{n} \partial_t \int_{\mathbb{R}^{2d}} v_{n-m}(t,q,p)\mu_m(q,p)\rho(q,p)dqdp = \sum_{m=0}^{n} \partial_t \int_{\mathbb{R}^{2d}} v_m(t,q,p)\mu_{n-m}(q,p)\rho(q,p)dqdp$$

$$= \sum_{m=0}^{n} \sum_{\ell=0}^{m} \int_{\mathbb{R}^{2d}} L_{m-\ell} v_\ell(t,q,p) \mu_{n-m}(q,p)\rho(q,p)dqdp$$

$$= \sum_{\ell=0}^{n} \int_{\mathbb{R}^{2d}} v_\ell(t,q,p) \sum_{m=0}^{n-\ell} L_m^* \mu_{n-\ell-m}(q,p)\rho(q,p)dqdp$$

$$= 0,$$

by definition of the coefficients $\mu_n$ (see (5.1)). Note that the computation above is justified because $\forall n$, $v_n$, $\mu_n$ and all their derivatives have polynomial growth. We deduce, for all $t \geq 0$,

$$\int_{\mathbb{R}^{2d}} \partial_t v_n(t,q,p)\rho(q,p)dqdp = -\sum_{m=1}^{n} \int_{\mathbb{R}^{2d}} \partial_t v_{n-m}(t,q,p)\mu_m(q,p)\rho(q,p)dqdp. \tag{5.3}$$

Now we compute the average of $F_n$. By (4.8), (4.9) and (5.3), we have for $t \geq 0$

$$\langle F_n(t) \rangle = \int_{\mathbb{R}^{2d}} F_n(t,q,p)\rho(q,p)dqdp = \int_{\mathbb{R}^{2d}} \partial_t v_n(t,q,p)\rho(q,p)dqdp - \int_{\mathbb{R}^{2d}} Lv_n(t,q,p)\rho(q,p)dqdp$$

$$= \int_{\mathbb{R}^{2d}} \partial_t v_n(t,q,p)\rho(q,p)dqdp$$

$$= -\sum_{m=1}^{n} \int_{\mathbb{R}^{2d}} \partial_t v_{n-m}(t,q,p)\mu_m(q,p)dqdp.$$

We rewrite (4.10) as follows: for all $(q,p) \in \mathbb{R}^{2d}$ and $t \geq 0$,

$$v_n(t,q,p) = \int_0^t \langle F_n(s) \rangle ds + \int_0^t P_{t-s}(F_n(s,q,p) - \langle F_n(s) \rangle) ds.$$

Using the previous expression obtained for $\langle F_n(s) \rangle$ and recalling the initial data for $v_n$, we deduce that for all $(q,p) \in \mathbb{R}^{2d}$ and $t \geq 0$

$$v_n(t,q,p) = -\sum_{m=1}^{n} \int_{\mathbb{R}^{2d}} v_{n-m}(t)\mu_m d\rho + \int \phi \mu_n d\rho$$

$$+ \int_0^t P_{t-s}(F_n(s,q,p) - \langle F_n(s) \rangle) ds.$$

Then, using $\int_{\mathbb{R}^{2d}} \mu_m(q,p)\rho(q,p)dqdp = 0$, for $m \in \mathbb{N}^*$ (Proposition 5.1), we get for all $(q,p) \in \mathbb{R}^{2d}$ and $t \geq 0$

$$v_n(t,q,p) - \int \phi \mu_n d\rho = -\sum_{m=1}^{n} \int_{\mathbb{R}^{2d}} \left( v_{n-m}(t,\tilde{q},\tilde{p}) - \int \phi \mu_{n-m} d\rho \right) \mu_m(\tilde{q},\tilde{p})\rho(\tilde{q},\tilde{p})d\tilde{q}d\tilde{p}$$
$$+ \int_0^t P_{t-s}(F_n(s,q,p) - \langle F_n(s) \rangle)ds.$$

Note that, since $L_\ell$, $\ell \in \mathbb{N}$ is a differential operator of order $2\ell + 2$ whose the coefficients belong to $\mathscr{C}_{pol}^\infty(\mathbb{R}^{2d})$ and contain no zero order terms, we have for $k \in \mathbb{N}$ that there exist $\gamma_{k,n} \in \mathbb{N}$ and $\zeta_{k,n,l} \in \mathbb{N}$ such that for $s \geq 0$,

$$\| F_n(s) - \langle F_n(s) \rangle \|_{k,\gamma_{k,n}} \leq \sum_{\ell=0}^{n-1} c_{k,\ell} |v_\ell(s)|_{b_{k,n,l},\zeta_{k,n,l}}$$
$$\leq \sum_{\ell=0}^{n-1} c_{k,\ell} \| v_\ell(s) - \int \phi \mu_\ell d\rho \|_{b_{k,n,l},\zeta_{k,n,\ell}},$$

where $b_{k,n,l} = k + 2(n-l) + 2$. We have used:

$$|v_n(t_{j+1},.)|_{\alpha,\beta} = |v_n(t_{j+1},.) - \int_{\mathbb{R}^{2d}} \phi(q,p)\rho_n(q,p)\rho(q,p)dqdp|_{\alpha,\beta}.$$

Then, using Proposition 2.15, we have for $k \in \mathbb{N}$ that there exists $r \in \mathbb{N}$ such that for $t \geq 0$

$$\| v_n(t) - \int \phi \mu_n d\rho \|_{k,r} \leq$$
$$\sum_{m=1}^{n} \left( \int_{\mathbb{R}^{2d}} |v_{n-m}(t,q,p) - \int \phi \mu_{n-m} d\rho|^2 \rho(q,p)dqdp \right)^{1/2} \left( \int_{\mathbb{R}^{2d}} |\mu_m(q,p)|^2 \rho(q,p)dqdp \right)^{1/2}$$
$$+ \int_0^t C_{k,n,i} e^{-\lambda(t-s)} \| F_n(s) - \langle F_n(s) \rangle \|_{k+(d+1),\gamma_{k+(d+1)}} ds.$$

Using the induction assumption, we have, for $k \in \mathbb{N}$ that there exists $r \in \mathbb{N}$ such that for $t \geq 0$,

$$\| v_n(t) - \int \phi \mu_n d\rho \|_{k,r} \leq \sum_{m=1}^{n} \tilde{c}_m P_{0,n-m,0}(t) e^{-\lambda t} \| \phi - \langle \phi \rangle \|_{(n+1)(d+1)+4n,\alpha_{0,j}}$$
$$+ C_{k,n} \sum_{l=0}^{n-1} c_{k,l} \int_0^t P_{k,l}(s) e^{-\lambda(t-s)} e^{-\lambda s} ds \| \phi - \langle \phi \rangle \|_{4n+k+(n+1)(d+1),\alpha_{k,j}}.$$

The conclusion follows.

The following Proposition ends the proof of our main result Theorem 2.17.

**Proposition 5.3** *Let $N$ be fixed and $\ell_N$ be fixed. Let $\delta_0 = \min(\frac{1}{\gamma}, \frac{\gamma\beta}{4\theta})$. Let $(q_k,p_k)$ be the discrete process defined by the implicit Euler scheme (2.12) or the implicit split-step scheme (2.11). Let $0 \leq \delta < \delta_0$, $\alpha_N = 6N + 2 + (N+1)(d+1)$ and $\phi \in \mathscr{C}_{pol}^\infty(\mathbb{R}^{2d}) \cap \mathscr{C}_{\ell_N}^{\alpha_N}(\mathbb{R}^{2d})$, then there exist strictly positive real number $C_N$ and an integer $k_N$ such that we have for $k \geq 0$,*

$$\| \mathbb{E}\phi(q_k,p_k) - v^{(N)}(t_k,.) \|_{0,k_N} \leq \delta^N C_N \| \phi - \langle \phi \rangle \|_{\alpha_N,\ell_N}, \tag{5.4}$$

where $t_k = k\delta$.

Moreover, for $0 \leq \delta < \delta_0$, $\lambda < \lambda_0$ and $\phi \in \mathscr{C}_{pol}^{\infty}(\mathbb{R}^{2d}) \cap \mathscr{C}_{\ell_N}^{\alpha_N}(\mathbb{R}^{2d})$, there exist a positive real number $C_N$, an integer $k_N$ and a positive polynomial function $P_N$ satisfying the following : For all $k \in \mathbb{N}$,

$$\| \mathbb{E}\phi(q_k, p_k) - \int \phi \mu^{(N)} d\rho \|_{0,k_N} \leq C_N \left( e^{-\lambda t_k} P_N(t_k) + \delta^N \right) \| \phi - \langle \phi \rangle \|_{\alpha_N, \ell_N},$$

where $t_k = k\delta$.

*Proof* Let $N$ be fixed and $l_N$ be fixed. Let $(q_k, p_k)$ be the discrete process defined by the implicit Euler scheme (2.12) or the implicit split-step scheme (2.11), $0 \leq \delta < \delta_0$ and $\phi \in \mathscr{C}_{pol}^{\infty}(\mathbb{R}^{2d}) \cap \mathscr{C}_{l_N}^{\alpha_N}(\mathbb{R}^{2d})$ be fixed. Let $(q,p) \in \mathbb{R}^{2d}$ such that $q_0 = q$ and $p_0 = p$. For all $k$, with $t_k = k\delta$, we have

$$\mathbb{E}\phi(q_k, p_k) - v^{(N)}(t_k, q, p) = \mathbb{E}v^{(N)}(0, q_k, p_k) - v^{(N)}(t_k, q, p)$$
$$= \mathbb{E} \sum_{j=0}^{k-1} \mathbb{E}^{q_{k-j-1}, p_{k-j-1}} \left( v^{(N)}(t_j, q_{k-j}, p_{k-j}) - v^{(N)}(t_{j+1}, q_{k-j-1}, p_{k-j-1}) \right).$$

Here we have used the notation $\mathbb{E}^{q_{k-j-1}, p_{k-j-1}}$ for the conditional expectation with respect to the filtration generated by $q_{k-j-1}$ and $p_{k-j-1}$. We have :

$$\mathbb{E}^{q_{k-j-1}, p_{k-j-1}} \left( v^{(N)}(t_j, q_{k-j}, p_{k-j}) - v^{(N)}(t_{j+1}, q_{k-j-1}, p_{k-j-1}) \right)$$
$$= \mathbb{E}^{q_{k-j-1}, p_{k-j-1}} \left( v^{(N)}(t_j, q_1(q_{k-j-1}), p_1(p_{k-j-1})) - v^{(N)}(t_{j+1}, q_{k-j-1}, p_{k-j-1}) \right),$$

where $(q_1(q), p_1(p))$ are the first step of the scheme (2.12) or of the scheme (2.11) when the initial condition is $(q,p)$. Using Theorem 4.1 with $t = t_j$, Proposition 2.9 and (5.2), we deduce that there exist integers $p_N$ and $k_N$ such that for all $k \in \mathbb{N}$

$$\| \mathbb{E}\phi(q_k, p_k) - v^{(N)}(t_k, .) \|_{0,k_N} \leq \delta^{N+1} C_N \sum_{j=0}^{k-1} \sup_{s \in ]0,\delta[, n=0,\ldots,N} |v_n(t_{j+1}, .)|_{2N+2, p_N}$$
$$\leq \delta^{N+1} C_N \| \phi - \langle \phi \rangle \|_{6N+2+(N+1)(d+1), \ell_N} \sum_{j=0}^{p-1} e^{\lambda t_j} P_N(t_j)$$
$$\leq \delta^{N+1} C_N \| \phi - \langle \phi \rangle \|_{6N+2+2(N+1)(d+1), \ell_N} \sum_{j=0}^{p-1} e^{-\tilde{\lambda} t_j},$$

for some constant $C_N$. We have used:

$$|v_n(t_{j+1}, .)|_{\alpha, \beta} = |v_n(t_{j+1}, .) - \int_{\mathbb{R}^{2d}} \phi(q, p) \rho_n(q, p) \rho(q, p) dq dp|_{\alpha, \beta}.$$

We conclude by using the fact that for a fixed constant $\gamma_1 > 0$, we have

$$\sum_{j=0}^{p-1} e^{-\gamma_1 j \delta} \leq \frac{1}{1 - e^{-\gamma_1 \delta}} \leq \frac{C}{\delta},$$

where the constant C depends on $\gamma_1$ and $\delta_0$. This shows (5.4). The second estimate is a consequence of (5.2).

# A Appendix: Proof of Proposition 2.15 and Lemma 2.16

## A.1 Notations, assumptions and result

In all this appendix, the constant $C$ may vary from line to line and we will use the following notations:
Let $\phi \in \mathscr{C}_{pol}^{\infty}(\mathbb{R}^{2d})$, then, for all $m \in \mathbb{N}$, there exists an integer $r_m$ such that $\phi \in \mathscr{C}_{r_m}^{m}(\mathbb{R}^{2d})$. For a multi-index $\mathbf{k} = (k_1,...,k_{2d}) \in \mathbb{N}^d$, we set $|\mathbf{k}| = k_1 + ... + k_{2d}$ and for a function $\phi \in C^{\infty}(\mathbb{R}^{2d})$, we set

$$D^{\mathbf{k}}\phi(x) = \frac{\partial^{|\mathbf{k}|}\phi(x)}{\partial_{x_1}^{k_1}...\partial_{x_{2d}}^{k_d}}, \quad x = (x_1,...,x_{2d}) \in \mathbb{R}^{2d}.$$

Moreover, we will assume in al the appendix that $\int \phi d\rho = 0$. We define $u$ by (2.26).

The aim of this appendix is to prove the following result (Proposition 2.15):

**Proposition A.1** *There exists a strictly positive real number $\lambda_0$ such that for any $m \in \mathbb{N}$, $\mathbf{k} \in \mathbb{N}^{2d}$ such that $|\mathbf{k}| = m$ and $\lambda < \lambda_0$, there exist an integer $s > 2r_{m+1+d}$ and a strictly positive real number $C$ such that for all $t \geq 0$ and $(q,p) \in \mathbb{R}^{2d}$*

$$\left|D^{\mathbf{k}}u(t,q,p)\right| \leq C(1+|q|^s+|p|^s) \| \phi \|_{m+d+1,r_{m+1+d}} \exp(-\lambda t).$$

To prove this Proposition, we will use the same idea as in [26]. Unlike in [26], we need to know how the estimate depend of $\phi$. Moreover, using an estimate of $u$ describes in [18], the proof is easier than in [26].
The proof proceeds as follows. We first show estimates on $u$ and its derivatives in an appropriate space. More precisely, we will show that for any $\lambda < 2\lambda_0$, $m \in \mathbb{N}$ and $\mathbf{k} \in \mathbb{N}^{2d}$ such that $|\mathbf{k}| = m$, there exists a strictly positive real number $C_m$ such that for all $t > 0$

$$\int |D^{\mathbf{k}}u(t,q,p)|^2 \pi_s(q,p) dq dp \leq C_m \exp(-\lambda t),$$

where the function $\pi_s$ is defined as

$$\pi_s = \frac{1}{\Gamma^s}, \tag{A.1}$$

for some integer $s$.
Moreover, we have the following result on $\pi_s$: For all multi-index $\mathbf{j}$ and integer $s$, there exists a function $\psi_{\mathbf{j},s} \in C^{\infty}$ such that

$$\partial^{\mathbf{j}}\pi_s(q,p) = \psi_{\mathbf{j},s}(q,p)\pi_s(q,p) \tag{A.2}$$

where

$$\psi_{\mathbf{j},s}(q,p) \xrightarrow{|(q,p)|\to\infty} 0.$$

Then, for any $m \in \mathbb{N}$ and $\mathbf{k} \in \mathbb{N}^{2d}$ such that $|\mathbf{k}| = m$, it is possible to choose an integer $s_m$ such that we have, for all $t > 0$, $s \geq s_m$ and $\zeta < 2\lambda_0$

$$\int_{\mathbb{R}^{2d}} |D^{\mathbf{k}}\big(u(t,q,p)\pi_s(q,p)\big)|^2 dq dp < C_m \exp(-\zeta t).$$

We then conclude by applying Sobolev imbedding Theorem (see [2]).

We will also use the following notation: for all $s \in \mathbb{N}$, $d\pi_s = \pi_s(q,p)dqdp$.

## A.2 Estimates on $u(t)$ and its derivatives in $L^2(\pi_s) = \{f : \mathbb{R}^{2d} \to \mathbb{R}; \int |f|^2 d\pi_s\}$

Using expression (2.24) of $L^{\top}$ and inequality (2.7) on $L\Gamma$, computations lead to

$$L^{\top}(\pi_s) = s\frac{L\Gamma}{\Gamma^{s+1}} - \frac{sd\sigma^2}{\Gamma}\pi_s + \frac{s(s+1)\sigma^2}{2}\frac{|\partial_p\Gamma|^2}{\Gamma^2}\pi_s + d\gamma\pi_s$$
$$\leq (-a_1 s + \gamma d)\pi_s + \Phi_s\pi_s,$$

where $\Phi_s(q,p) \xrightarrow{|(q,p)|\to+\infty} 0$ and $a_1$ is defined by (2.7). Hence, for each $s \in \mathbb{N}^*$, there exists a real number $\nu_s > 0$ such that

$$L^{\top}(\pi_s) \leq \nu_s \pi_s. \tag{A.3}$$

We will now show the following proposition then we will use Sobolev inequalities to prove Proposition A.1:

**Proposition A.2** *Let $m \in \mathbb{N}$, $\zeta < 2\lambda$ and $\mathbf{k} \in \mathbb{N}^{2d}$ such that $|\mathbf{k}| = m$ be fixed. There exists an integer $s_m > 2r_m$ such that, for all $s \geq s_m$, there exists a strictly positive real number $C_{m,s}$ such that we have for all $T > 0$,*

$$\int |D^{\mathbf{k}} u(T)|^2 d\pi_s \leq C_{m,s} \| \phi \|^2_{m,r_m} \exp(-\zeta T), \tag{A.4}$$

*where $r_m$ is defined at the beginning of this appendix.*

Proposition A.2 is a corollary of the following result:

**Proposition A.3** *For all $m \in \mathbb{N}$, $\mathbf{k} \in \mathbb{N}^{2d}$ such that $|\mathbf{k}| = m$, positive polynomial function $Q$ and $\zeta < 2\lambda$, there exists an integer $s_m > 2r_m$ such that for all $s \geq s_m$ there exists a strictly positive real numbers $C_{m,Q,s}$ such that we have for all $T > 0$*

$$\int_0^T \exp(\zeta t) \int Q |D^{\mathbf{k}} u(t)|^2 d\pi_s dt \leq C_{m,Q,s} \| \phi \|^2_{m,r_m}. \tag{A.5}$$

For the convenience of the reading, we will first show Propositions A.2 and A.3 for $m = 0$. Then, we will show Proposition A.3 for $m = 1$ and explain how to deduce Proposition A.2 for $m = 1$. Finally, we will show Proposition A.3 by induction and show that Proposition A.2 is a corollary of Proposition A.3. The idea to prove Proposition A.3 for $m \geq 1$ and $\mathbf{k} \in \mathbb{N}^{2d}$ such that $|\mathbf{k}| = m$ is the following: We first show the result for $Q|\partial_p D^{\mathbf{k-1}}.|^2$ and $Q|(\alpha \partial_p - \partial_q) D^{\mathbf{k-1}}.|^2$. We can then deduce the result for $Q|\partial_q D^{\mathbf{k-1}}.|^2$.

To show this two Propositions for $m = 0$, we need of a better point-wise estimatie of $u$:

**Lemma A.4** *There exists $C = C(r_0) > 0$, $\lambda_0 = \lambda_0(r_0) > 0$ such that, for all $t \geq 0$ and $(q,p) \in \mathbb{R}^{2d}$,*

$$|u(t,q,p)| \leq C \Gamma^{r_0}(q,p) \exp(-\lambda_0 t) \| \phi \|_{0,r_0}. \tag{A.6}$$

A proof of this result can be found in [18] or [10]. To show this Lemma, the two main ingredients are the property (2.6) on $\Gamma$ and that for $x \in \mathbb{R}^{2d}$, $t > 0$ and open $\mathscr{O} \subset \mathbb{R}^{2d}$, the transition kernel for (2.4) satisfies $Q_t(x, \mathscr{O}) > 0$. Under Assumption **B**, the second property is true (see [18]).

We have the following Corollary:

**Lemma A.5** *For $s > 2r_0$, there exists a strictly positive real number $C_s$ such that*

$$\int |u(t)|^2 d\pi_s \leq C_s \| \phi \|^2_{0,r_0} \exp(-2\lambda t) \quad \forall t \geq 0. \tag{A.7}$$

*Moreover, for $Q$ a positive polynomial function, we have, for all $s > s_Q \geq 2r_0$, that there exists a strictly positive real number $C_{Q,s}$ such that*

$$\int Q|u(t)|^2 d\pi_s \leq C_{Q,s} \| \phi \|^2_{0,r_0} \exp(-2\lambda t) \quad \forall t \geq 0. \tag{A.8}$$

Then, the results (A.4) and (A.5) for $m = 0$ are a consequence of Lemma A.5.

The following Lemma is the key of the proof of all the other Lemmas.

**Lemma A.6** *Let $A$ be a linear operator and $Q$ a polynomial function. There exists an integer $s_{A,Q}$ such that for all $s \geq s_{A,Q}$, we have for all $\zeta > 0$ and $T > 0$*

$$\exp(\zeta T) \int Q |Au(t)|^2 d\pi_s + \sigma^2 \int_0^T \exp(\zeta t) \int Q |\partial_p (Au(t))|^2 d\pi_s dt$$
$$\leq \int Q |Au(0)|^2 d\pi_s + (\zeta + \nu_s) \int_0^T \exp(\zeta t) \int Q |Au(t)|^2 d\pi_s dt$$
$$+ 2 \int_0^T \exp(\zeta t) \int Q \langle [A,L]u(t), Au(t) \rangle d\pi_s dt - \int_0^T \exp(\zeta t) \int |Au(t)|^2 LQ d\pi_s dt$$
$$- \sigma^2 \int_0^T \exp(\zeta t) \int \langle \partial_p Q, \partial_p |Au(t)|^2 \rangle d\pi_s dt, \tag{A.9}$$

*where, for $A$ and $B$ two linear operators, $[A,B] = AB - BA$.*

*Proof* Let *s* large enough such that $\int Q|Au(0)|^2 d\pi_s < \infty$. Using (2.23) and (2.25), we get

$$\begin{aligned}\frac{d}{dt}[\exp(\zeta t)Q|Au(t)|^2] &= \zeta \exp(\zeta t)Q|Au(t)|^2 + 2\exp(\zeta t)Q\langle ALu(t), Au(t)\rangle \\ &= \zeta\exp(\zeta t)Q|Au(t)|^2 + 2\exp(\zeta t)Q\langle LAu(t), Au(t)\rangle \\ &\quad + 2\exp(\zeta t)Q\langle [A,L]u(t), Au(t)\rangle \\ &= \zeta\exp(\zeta t)Q|Au(t)|^2 + \exp(\zeta t)QL|Au(t)|^2 - \sigma^2\exp(\zeta t)Q|\partial_p(Au(t))|^2 \\ &\quad + 2\exp(\zeta t)Q\langle [A,L]u(t), Au(t)\rangle \\ &= \zeta\exp(\zeta t)Q|Au(t)|^2 + \exp(\zeta t)L(Q|Au(t)|^2) - \exp(\zeta t)|Au(t)|^2 LQ \\ &\quad - \sigma^2\exp(\zeta t)\langle \partial_p Q, \partial_p|Au(t)|^2\rangle - \sigma^2\exp(\zeta t)Q|\partial_p(Au(t))|^2 \\ &\quad + 2\exp(\zeta t)Q\langle [A,L]u(t), Au(t)\rangle.\end{aligned}$$

We integrate with respect to *t*,

$$\begin{aligned}\exp(\zeta T)Q|Au(T)|^2 &= Q|Au(0)|^2 + \zeta\int_0^T \exp(\zeta t)Q|Au(t)|^2 dt + \int_0^T \exp(\zeta t)L(Q|Au(t)|^2)dt \\ &\quad - \int_0^T \exp(\zeta t)|Au(t)|^2 LQ dt - \sigma^2 \int_0^T \exp(\zeta t)\langle \partial_p Q, \partial_p|Au(t)|^2\rangle dt \\ &\quad - \sigma^2\int_0^T \exp(\zeta t)Q|\partial_p Au(t)|^2 dt + 2\int_0^T \exp(\zeta t)Q\langle [A,L]u(t), Au(t)\rangle dt.\end{aligned}$$

We integrate with respect to $\pi_s$. Using the inequality (A.3) on $L^\top \pi_s$, we have (A.13).

We will also need of following computations:

**Lemma A.7** *Let $k \in \mathbb{N}$, we have for any $(q,p) \in \mathbb{R}^{2d}$*

$$-LH^k(q,p) \leq 2\gamma k H^k(q,p) \tag{A.10}$$

*and for all $t \geq 0$ and linear operator A,*

$$\langle \partial_p H^k, \partial_p|Au(t)|^2\rangle \leq 2kH^k|Au(t)|^2 + 2kH^{k-1}|\partial_p Au(t)|^2. \tag{A.11}$$

*Proof* For any $(q,p) \in \mathbb{R}^{2d}$ and $k \in \mathbb{N}^*$, we have

$$\partial_p H^k(q,p) = kH^{k-1}(q,p)p$$

and

$$-LH^k(q,p) = \gamma k|p|^2 H^{k-1}(q,p) - \frac{k}{2}H^{k-1}(q,p) - \frac{k(k-1)}{2}|p|^2 H^{k-2}(q,p)$$
$$\leq \gamma k|p|^2 H^{k-1}(q,p) \leq 2\gamma k H^k(q,p).$$

We have used the positivity of *V*. Moreover, we have for any $t \geq 0$, $\ell \in \mathbb{N}$ and for each component of *Au* that we still write *Au*,

$$\begin{aligned}\langle \partial_p H^{2l}, \partial_p|Au(t)|^2\rangle &= 4\ell H^\ell Au(t)\langle pH^{\ell-1}, \partial_p Au(t)\rangle \\ &\leq 2\ell H^{2\ell}|Au(t)|^2 + 2\ell|p|^2 H^{2\ell-2}|\partial_p Au(t)|^2 \\ &\leq 2\ell H^{2\ell}|Au(t)|^2 + 4\ell H^{2\ell-1}|\partial_p Au(t)|^2\end{aligned}$$

and

$$\begin{aligned}\langle \partial_p H^{2l+1}, \partial_p|Au(t)|^2\rangle &= 2(2\ell+1)\langle H^\ell u(t)p, H^\ell \partial_p Au(t)\rangle \\ &\leq (2\ell+1)H^{2\ell}|p|^2|u(t)|^2 + (2\ell+1)H^{2\ell}|\partial_p Au(t)|^2 \\ &\leq 2(2\ell+1)H^{2\ell+1}|u(t)|^2 + (2\ell+1)H^{2\ell}|\partial_p Au(t)|^2.\end{aligned}$$

We now show the results (A.4) and (A.5) for $m=1$. First, we show the two following preliminary Lemmas.

**Lemma A.8** *Let $Q$ be a positive polynomial function. Let $s > 2r_0$ large enough and $\zeta < 2\lambda_0$. There exists a strictly positive real number $C_{Q,s}$ such that for all $T > 0$*

$$\int_0^T \exp(\zeta t) \int Q |\partial_p u(t)|^2 d\pi_s dt \leq C_{Q,s} \| \phi \|_{0,r_0}^2 . \tag{A.12}$$

*Proof* Let $\zeta > 0$. We have $[Id, L] = 0$, then, using Lemma A.6 with $A = I_d$ and $Q = 1$, we get for $T > 0$

$$\exp(\zeta T) \int |u(t)|^2 d\pi_s + \sigma^2 \int_0^T \exp(\zeta t) \int |\partial_p u(t)|^2 d\pi_s dt$$
$$\leq \int |u(0)|^2 d\pi_s + (\zeta + v_s) \int_0^T \exp(\zeta t) \int |u(t)|^2 d\pi_s dt.$$

We choose $\zeta < 2\lambda_0$ and, using (A.7), we bound the last term. Then, we have (A.12) for $Q = 1$.

Let $Q$ a positive polynomial function. Using Lemma A.6 for $A = I_d$, we get for $s$ large enough, $\zeta > 0$ and $T > 0$

$$\exp(\zeta T) \int Q |u(t)|^2 d\pi_s + \sigma^2 \int_0^T \exp(\zeta t) \int Q |\partial_p (u(t))|^2 d\pi_s dt$$
$$\leq \int Q |u(0)|^2 d\pi_s + (\zeta + v_s) \int_0^T \exp(\zeta t) \int Q |u(t)|^2 d\pi_s dt$$
$$- \int_0^T \exp(\zeta t) \int |u(t)|^2 L Q d\pi_s dt - \sigma^2 \int_0^T \exp(\zeta t) \int \langle \partial_p Q, \partial_p |u(t)|^2 \rangle d\pi_s dt, \tag{A.13}$$

As $V \in \mathscr{C}_{pol}^\infty(\mathbb{R}^d)$, there exists a positive polynomial function $Q_1$ such that $|LQ| \leq Q_1$, then we have for any $T > 0$ and $s$ large enough

$$- \int_0^T \exp(\zeta t) \int |u(t)|^2 L Q d\pi_s dt - \sigma^2 \int_0^T \exp(\zeta t) \int \langle \partial_p Q, \partial_p |u(t)|^2 \rangle d\pi_s dt$$
$$\leq \int_0^T \exp(\zeta t) \left(2 Q_1 + \sigma^2 |\partial_p Q|^2\right) |u(t)|^2 d\pi_s dt + \sigma^2 \int_0^T \exp(\zeta t) \int |\partial_p u(t)|^2 d\pi_s dt.$$

We choose $0 < \zeta < 2\lambda_0$ and use (A.8) and (A.12) with $Q = 1$ to have the result (A.12) for any positive polynomial function $Q$.

**Lemma A.9** *Let $s > 2r_1$ large enough and $\zeta < 2\lambda_0$. There exists a strictly positive real number $\alpha_s$ such that for $\alpha > \alpha_s$, there exists a positive real number $C_s$ such that, for all $T > 0$,*

$$\int_0^T \exp(\zeta t) \int |\alpha \partial_p u(t) - \partial_q u(t)|^2 d\pi_s dt + \int_0^T \exp(\zeta t) \int |\partial_p (\alpha \partial_p u(t) - \partial_q u(t))|^2 d\pi_s dt$$
$$\leq C_s \| \phi \|_{1,r_1}^2 . \tag{A.14}$$

*Let $k > 0$ an integer. For any $\zeta < 2\lambda_0$ and $s > 2r_1$ large enough, there exists a strictly positive real number $\alpha_k$ such that $\forall \alpha \geq \alpha_k$ there exists a strictly positive real number $C_{k,s}$ such that, for all $T > 0$,*

$$\exp(\zeta T) \int H^k |\alpha \partial_p u(T) - \partial_q u(T)|^2 d\pi_s + \int_0^T \exp(\zeta t) \int H^k |\partial_p (\alpha \partial_p u(t) - \partial_q u(t))|^2 d\pi_s dt$$
$$\leq C_{k,s} \| \phi \|_{1,r_1}^2 . \tag{A.15}$$

*Proof* Let $\alpha > 0$ and $i \in \{1, ..., d\}$. We use the following notation: $A_1 = \alpha \partial_{p_i} - \partial_{q_i}$. We have for any function $\psi \in C^\infty(\mathbb{R}^{2d})$

$$[A_1, L]\psi = -\alpha A_1 \psi + \alpha(\alpha - \gamma) \partial_{p_i} \psi + \langle \partial_{q_i} \partial_q V, \partial_p \psi \rangle.$$

Using the polynomial growth of $\partial_{q_i} \partial_q V$, we have that there exist a positive polynomial function $Q_1$, $\varepsilon_1 > 0$ and $\varepsilon_2 > 0$ such that for any function $\psi \in C^\infty(\mathbb{R}^{2d})$

$$2\langle [A_1, L]\psi, A_1 \psi \rangle = -2\alpha |A_1 \psi|^2 + 2 A_1 \psi \langle \partial_{q_i} \partial_q V, \partial_p \psi \rangle + 2\alpha(\alpha - \gamma) \partial_{p_i} \psi A_1 \psi$$
$$\leq (\varepsilon_1 + \varepsilon_2 - 2\alpha) |A_1 \psi|^2 + \frac{Q_1}{\varepsilon_1} |\partial_p \psi|^2 + \frac{\alpha(\alpha - \gamma)^2}{\varepsilon_2} |\partial_{p_i} \psi|^2. \tag{A.16}$$

Then, choosing $s$ large enough and using Lemma A.6 with $Q = 1$, we get for all $T > 0$, $\varepsilon_1$, $\varepsilon_2$ and $\zeta > 0$

$$\exp(\zeta T) \int |A_1 u(T)|^2 d\pi_s + \sigma^2 \int_0^T \exp(\zeta t) \int |\partial_p(A_1 u(t))|^2 d\pi_s dt$$
$$\leq \int |A_1 u(0)|^2 d\pi_s + \frac{1}{\varepsilon_1} \int_0^T \exp(\zeta t) \int Q_1 |\partial_p u(t)|^2 d\pi_s dt$$
$$+ \frac{(\alpha(\alpha-\gamma))^2}{\varepsilon_2} \int_0^T \exp(\zeta t) \int |\partial_{p_i} u(t)|^2 d\pi_s dt$$
$$+ (\nu_s + \zeta + \varepsilon_1 + \varepsilon_2 - 2\alpha) \int_0^T \exp(\zeta t) \int |A_1 u(t)|^2 d\pi_s dt.$$

We choose $\zeta < 2\lambda_0$ and $\alpha$, $\varepsilon_1$ and $\varepsilon_2$ such that $\nu_s + \zeta + \varepsilon_1 + \varepsilon_2 - 2\alpha < 0$. We then use (A.12) to prove (A.14).

We now prove (A.15) by recursion on $k$. We have proved the case $k = 0$. Let us assume (A.15) is true for $k - 1$. We want to obtain it for $k$.

Using computations for $k = 0$, (A.16), (A.10), (A.11) and Lemma A.6, we get for $s$ large enough, $\varepsilon_1 > 0$, $\varepsilon_2 > 0$, $\zeta > 0$ and $T > 0$

$$\exp(\zeta T) \int H^k |\alpha \partial_p u(t) - \partial_q u(t)|^2 d\pi_s + \sigma^2 \int_0^T \exp(\zeta t) \int H^k |\partial_p(\alpha \partial_p u(t) - \partial_q u(t))|^2 d\pi_s dt$$
$$\leq \int H^k |\alpha \partial_p u(0) - \partial_q u(0)|^2 d\pi_s$$
$$+ (\zeta + \varepsilon_2 + \varepsilon_1 - 2\alpha + 2k\sigma^2 + 2k\gamma + \nu_s) \int_0^T \exp(\zeta t) \int H^k |\alpha \partial_p u(t) - \partial_q u(t)|^2 d\pi_s dt$$
$$+ 2k\sigma^2 \int_0^T \exp(\zeta t) \int H^{k-1} |\partial_p(\partial_p u(t) - \partial_q u(t))|^2 d\pi_s dt$$
$$+ \left(\frac{Q_1}{\varepsilon_1} + \frac{\alpha^2(\alpha+\gamma)^2}{\varepsilon_2}\right) \int_0^T \exp(\zeta t) \int H^k |\partial_p u(t)|^2 d\pi_s dt,$$

where $Q_1$ is a positive polynomial function. We take $\zeta < 2\lambda$ and choose $\alpha$, $\varepsilon_2$ and $\varepsilon_1$ such that $\varepsilon_2 + \varepsilon_1 - 2\alpha + 2k\sigma^2 + 2k\gamma + \zeta + \nu_s \leq 0$. Then, using the induction hypothesis on $k$, the polynomial growth of H and (A.12), we obtain (A.15) for $k$.

**Remark A.10** *Using the fact that $q^2 \leq CH(q,p)$ and $p^2 \leq 2H(q,p)$, (A.14) and (A.15), we have for $\zeta < 2\lambda_0$, $s$ large enough and $Q$ a positive polynomial function that there exist real positive numbers $\alpha_s$ and $C_{Q,s}$ depending also of $\alpha_s$, such that, for $T > 0$,*

$$\int_0^T \exp(\zeta t) \int Q |\alpha \partial_p u(t) - \partial_q u(t)|^2 d\pi_s dt \leq C_{Q,s} \| \phi \|_{1,r_1}^2 . \tag{A.17}$$

Using above Remark, we can show the following Lemma:

**Lemma A.11** *Let $Q$ be a positive polynomial function. Let $s > 2r_1$ large enough. For all $\zeta < 2\lambda_0$, there exists a strictly positive real number $C_{Q,s}$ such that for all $T > 0$*

$$\int_0^T \exp(\zeta t) \int Q |\partial_q u(t)|^2 d\pi_s dt \leq C_{Q,s} \| \phi \|_{1,r_1}^2 . \tag{A.18}$$

*Proof* Let $\zeta < 2\lambda_0$. We have the following inequality: for any function $\psi \in C^\infty(\mathbb{R}^{2d})$ and $\alpha > 0$

$$|\partial_q \psi|^2 \leq |(\partial_q - \alpha \partial_p)\psi|^2 + \alpha^2 |\partial_p \psi|^2.$$

Then, using (A.17) and (A.12), we have (A.18).

Using (A.12) and (A.18), we obtain (A.5) for $m = 1$. We will now show (A.4) for $m = 1$.

**Lemma A.12** *For $s > 2r_1$ large enough and $\zeta < 2\lambda_0$, there exists a strictly positive real number $C_s$ such that for $T > 0$,*

$$\int |\partial_p u(T)|^2 d\pi_s \leq C_s \| \phi \|_{1,r_1}^2 \exp(-\zeta T).$$

*Proof* We have for $t \geq 0$

$$2\langle [\partial_{p_i},L]u(t), \partial_{p_i}u(t)\rangle = 2\partial_{q_i}u(t)\partial_{p_i}u(t) - 2\gamma|\partial_{p_i}u(t)|^2 \leq |\partial_{q_i}u(t)|^2 + (1-2\gamma)|\partial_{p_i}u(t)|^2.$$

Then, using Lemma A.6, we get for $T > 0$ and $\zeta > 0$

$$\exp(\zeta t)\int |\partial_p u(t)|^2 d\pi_s \leq \int |\partial_p u(0)|^2 d\pi_s + (\zeta + \nu_s + 1 - 2\gamma)\int_0^T \exp(\zeta t)\int |\partial_p u(t)|^2 d\pi_s ds$$
$$+ \int_0^T \exp(\zeta t)\int |\partial_q u(t)|^2 d\pi_s ds. \quad (A.19)$$

We choose $\zeta < 2\lambda_0$ and use (A.12) and (A.18) to conclude.

**Lemma A.13** *For $s > 2r_1$ large enough and $\zeta < 2\lambda_0$, there exists a strictly positive real number $C_s$ such that for all $T > 0$,*

$$\int |\partial_q u(T)|^2 d\pi_s \leq C_s \|\phi\|_{1,r_1}^2 \exp(-\zeta T).$$

*Proof* We have for $t \geq 0$

$$2\langle [\partial_{q_i},L]u(t), \partial_{q_i}u(t)\rangle = 2\langle \partial_{q_i}DV, \partial_p u(t)\rangle \partial_{q_i}u(t) \leq Q|\partial_p u(t)|^2 + |\partial_{q_i}u(t)|^2,$$

where $Q$ is a positive polynomial function such that $|\partial_{q_i}DV|^2 < Q$. Then, we use Lemma A.6, $\zeta < 2\lambda_0$, (A.12) and (A.18) to conclude.

We have shown the result (A.4) for $m = 1$.

We will now show equation (A.5) by induction on $m$. We already proved the case $m = 1$. We suppose that (A.5) holds up to $m$ and we want to obtain it for $m+1$. First, we show a result on $\partial_p D^{\mathbf{k}} u(t)$ where $\mathbf{k} \in \mathbb{N}^{2d}$ such that $|\mathbf{k}| = m$.

**Lemma A.14** *Let us assume that induction hypothesis (A.5) holds up to $m$. For $s > 2r_m$ large enough, $\zeta < 2\lambda_0$ and $\mathbf{k} \in \mathbb{N}^{2d}$ such that $|\mathbf{k}| = m$, there exists a strictly positive real number $C_{m,s}$ such that for all $T > 0$,*

$$\int_0^T \exp(\zeta t)\int |\partial_p D^{\mathbf{k}} u(t)|^2 d\pi_s dt \leq C_{m,s} \|\phi\|_{m,r_m}^2. \quad (A.20)$$

*Proof* By induction, we can show for any function $\psi$

$$[D^{\mathbf{k}}, L]\psi \leq \sum_{\mathbf{i} \in \mathbb{N}^{2d}, |\mathbf{i}| \leq m} P_{\mathbf{i}} |D^{\mathbf{i}}\psi|,$$

where $P_{\mathbf{i}}$ is a positive polynomial function which depend of the polynomial growth of $V$ and its derivatives. Then, we get for $t > 0$

$$2\langle [D^{\mathbf{k}}, L]u(t), D^{\mathbf{k}} u(t)\rangle \leq \sum_{\mathbf{i} \in \mathbb{N}^{2d}, |\mathbf{i}| \leq m} P_{\mathbf{i}}^2 |D^{\mathbf{i}} u(t)|^2 + C|D^{\mathbf{k}} u(t)|^2, \quad (A.21)$$

where $C$ is a positive real number. Using Lemma A.6, we get for $T > 0$ and $\zeta > 0$

$$\exp(\zeta T)\int |D^{\mathbf{k}} u(t)|^2 d\pi_s + \sigma^2 \int_0^T \exp(\zeta t)\int |\partial_p(D^{\mathbf{k}} u(t))|^2 d\pi_s dt$$
$$\leq \int |D^{\mathbf{k}} u(0)|^2 d\pi_s + (\zeta + \nu_s + C)\int_0^T \exp(\zeta t)\int |D^{\mathbf{k}} u(t)|^2 d\pi_s dt$$
$$+ \sum_{\mathbf{i} \in \mathbb{N}^{2d}, |\mathbf{i}| \leq m} \int_0^T \exp(\zeta t)\int P_{\mathbf{i}}^2 |D^{\mathbf{i}} u(t)| d\pi_s dt,$$

We choose $\zeta < 2\lambda_0$ and use induction hypothesis (A.5) to conclude.

**Lemma A.15** *Let us assume that induction hypothesis (A.5) holds. For $\zeta < 2\lambda_0$, $s > 2r_m$ large enough, $\ell \in \mathbb{N}$ and $\mathbf{k} \in \mathbb{N}^{2d}$ such that $|\mathbf{k}| = m$, there exists a strictly positive real number $C_{m,\ell,s}$ such that*

$$\int_0^T \exp(\zeta t)\int H^\ell |\partial_p D^{\mathbf{k}} u(t)|^2 d\pi_s dt \leq C_{m,\ell,s} \|\phi\|_{m,r_m}^2, \quad (A.22)$$

*for all $T > 0$.*

*Proof* We proceed by induction on $\ell$.

We also have the result for $\ell = 0$ (see (A.20)). Suppose that the induction hypothesis (A.22) holds for $\ell - 1$. Using computations for $\ell = 0$ (A.21), (A.10), (A.11) and Lemma A.6, we get for $s$ large enough, $\varepsilon_1 > 0$, $\varepsilon_2 > 0$, $\zeta > 0$ and $T > 0$

$$\exp(\zeta T) \int H^\ell |D^{\mathbf{k}} u(t)|^2 d\pi_s + \sigma^2 \int_0^T \exp(\zeta t) \int H^\ell |\partial_p D^{\mathbf{k}} u(t)|^2 d\pi_s dt$$

$$\leq \int H^\ell |D^{\mathbf{k}} u(0)|^2 d\pi_s$$

$$+ (\zeta + C + 2\ell\sigma^2 + 2\ell\gamma + v_s) \int_0^T \exp(\zeta t) \int H^\ell |D^{\mathbf{k}} u(t)|^2 d\pi_s dt$$

$$+ 2\ell\sigma^2 \int_0^T \exp(\zeta t) \int H^{\ell-1} |\partial_p D^{\mathbf{k}} u(t)|^2 d\pi_s dt$$

$$+ \sum_{\mathbf{i} \in \mathbb{N}^{2d}, |\mathbf{i}| \leq m} \int_0^T \exp(\zeta t) \int H^\ell P_{\mathbf{i}}^2 |D^{\mathbf{i}} u(t)| d\pi_s dt,$$

where $C$ and $P_{\mathbf{i}}$ are defined in (A.21). We take $\zeta < 2\lambda_0$. Then, using the induction hypothesis on $\ell$, the polynomial growth of H and induction hypothesis on $m$, we obtain (A.15) for $\ell$.

**Remark A.16** *Using the same ideas as in Remark A.10, we prove the following result. Let $Q$ a positive polynomial function. For $s > 2r_m$ large enough, $\zeta < 2\lambda_0$ and $\mathbf{k} \in \mathbb{N}^{2d}$ such that $|\mathbf{k}| = m$, there exists a strictly positive real number $C_{m,Q,s}$ such that for all $T > 0$*

$$\exp(\zeta T) \int Q |D^{\mathbf{k}} u(T)|^2 d\pi_s + \int_0^T \exp(\zeta t) \int Q |\partial_p D^{\mathbf{k}} u(t)|^2 d\pi_s dt \leq C_{m,Q,s} \| \phi \|_{m,r_m}^2. \tag{A.23}$$

**Lemma A.17** *Let us assume that induction hypothesis (A.5) holds. Let $D_q^m u(t)$ denote an arbitrary partial derivative of $u(t)$ of the type $\partial_{q_{i_1}}...\partial_{q_{i_m}}$. For $s > 2r_{m+1}$ large enough and $\zeta < 2\lambda_0$, there exists a strictly positive real number $\tilde{\alpha}_0$ such that for all $\alpha \geq \tilde{\alpha}_0$, there exists a strictly positive real number $C$ depending of $m$, $s$ and $\alpha$ such that for all $T > 0$,*

$$\int_0^T \exp(\zeta t) \int |\alpha \partial_p (D_q^m u(t)) - \partial_q (D_q^m u(t))|^2 d\pi_s dt \leq C \| \phi \|_{m+1,r_{m+1}}^2. \tag{A.24}$$

*Proof* Let $\alpha > 0$, $i \in \{1,...,d\}$ and $A_1 = \alpha \partial_p - \partial_q$. For any function $\psi \in C^\infty(\mathbb{R}^{2d})$, we have

$$[A_1 D_q^m, L] \psi = A_1 [D_q^m, L] \psi + [A_1, L](D_q^m \psi).$$

Moreover, by induction an $m$, we can show for any function $\psi \in C^\infty(\mathbb{R}^{2d})$

$$A_1 [D_q^m, L] \psi \leq \sum_{\mathbf{i} \in \mathbb{N}^{2d}, |\mathbf{i}| \leq m} P_{\mathbf{i}} |\partial_p D^{\mathbf{i}} \psi|,$$

where $P_{\mathbf{i}}$ is a positive polynomial function which depend of the polynomial growth of $V$ and its derivatives. Then, we get for $t > 0$

$$2 \langle A_1 [D_q^m, L] u(t), A_1 D_q^m u(t) \rangle \leq \sum_{\mathbf{i} \in \mathbb{N}^{2d}, |\mathbf{i}| \leq m} P_{\mathbf{i}}^2 |\partial_p D^{\mathbf{i}} u(t)|^2 + C |A_1 D_q^m u(t)|^2, \tag{A.25}$$

where $C$ is a positive real number. Using (A.25) and (A.16), we have that there exist a positive polynomial function $Q_1$, $\varepsilon_1 > 0$ and $\varepsilon_2 > 0$ such that

$$\langle [A_1 D_q^m, L] u(t), A_1 D_q^m u(t) \rangle \leq \sum_{\mathbf{i} \in \mathbb{N}^{2d}, |\mathbf{i}| \leq m} P_{\mathbf{i}}^2 |\partial_p D^{\mathbf{i}} u(t)|^2 + (\varepsilon_1 + \varepsilon_2 - 2\alpha + C) |A_1 D_q^m u(t)|^2$$

$$+ \frac{Q_1}{\varepsilon_1} |\partial_p D_q^m u(t)|^2 + \frac{\alpha(\alpha - \gamma)^2}{\varepsilon_2} |\partial_{p_i} D_q^m u(t)|^2.$$

To conclude, we proceed as Lemma A.9: We choose $s$ large enough and use Lemma A.6. We the choose $\zeta < 2\lambda$ and $\alpha$, $\varepsilon_1$ and $\varepsilon_2$ such that $v_s + \zeta + \varepsilon_1 + \varepsilon_2 - 2\alpha + C < 0$ and use (A.23).

**Lemma A.18** *Let us assume that induction hypothesis (A.5) holds. Let $Q$ a positive polynomial function. For $s > 2r_{m+1}$ large enough and $\zeta < 2\lambda_0$, there exists a strictly positive real number $\tilde{\alpha}_s$ such that for all $\alpha \geq \tilde{\alpha}_s$, there exists a strictly positive real number $C_{m,s}$ which also depends of $\alpha$ and $\zeta$ such that for all $T > 0$,*

$$\int_0^T \exp(\zeta t) \int Q |\alpha \partial_p (D_q^m u(t)) - \partial_q (D_q^m u(t))|^2 d\pi_s dt \leq C_{m,s} \| \phi \|_{m+1,r_{m+1}}^2. \tag{A.26}$$

*Proof* First, we prove this result by induction on $k$ for $H^k$. We use same arguments and computations used to prove (A.15) and (A.22), then we have the result for $Q$ (see Remark A.10).

**Remark A.19** Using (A.23) and (A.26), we prove that inequality (A.5) holds up to $m+1$. Indeed, if there is a derivative in a direction $p_1,...,p_d$, then we use the result (A.23). In the other case, we use (A.26), (A.23) and $\partial_q D_q^m = -(\alpha \partial_{p_i} D_q^m - \partial_{q_i} D_q^m) + \alpha \partial_{p_i} D_q^m$, for $i \in \{1,...,d\}$.

We have shown that (A.5) is true for all $m$, then we have shown Proposition A.3. We now prove that (A.4) is true for the derivatives of order $m$, where $m \in \mathbb{N}^*$.

**Lemma A.20** *Let $m \in \mathbb{N}^*$ be fixed. For $s > 2r_m$ large enough, $\zeta < 2\lambda_0$ and $\mathbf{k} \in \mathbb{N}^{2d}$ such that $|\mathbf{k}| = m$, there exists a strictly positive real number $C_{m,s}$ such that for all $T > 0$,*

$$\int |D^{\mathbf{k}}u(T)|^2 d\pi_s \leq C_{m,s} \| \phi \|_{m,r_m}^2 \exp(-\zeta T).$$

*Proof* We proceed as to prove Lemma A.11 and Lemma A.12 and we use Lemma A.6, (A.21) and Proposition A.3.

We have shown the result (A.4).

## A.3 Point-wise estimates of $u$ and its derivatives

*Proof (Proof of Proposition A.1)* Using (A.4) and (A.2), we get : For all $m \in \mathbb{N}$, $n \in \mathbb{N}$, $\mathbf{k} \in \mathbb{N}^{2d}$ such that $|\mathbf{k}| = m$, $\mathbf{l} \in \mathbb{N}^{2d}$ such that $|\mathbf{l}| = m$, $\zeta < 2\lambda_0$ and $s$ large enough, there exists a strictly positive real number $C_{m,n,s}$ such that for all $t > 0$

$$\int \left| D^{\mathbf{l}}\left(D^{\mathbf{k}}u(t,q,p)\pi_s(q,p)\right) \right|^2 dqdp \leq C_{m,n,s} \| \phi \|_{m+n,r_{m+n}}^2 \exp(-\zeta t).$$

Then, using Sobolev embedding theorem [2], we get, for $n = d+1$

$$|D^{\mathbf{k}}u(t,q,p)|^2 \pi_s^2(q,p) \leq C_{m,d+1,s} \| \phi \|_{m+d+1,r_{m+d+1}}^2 \exp(-\zeta t),$$

for all $t > 0$ and $(q,p) \in \mathbb{R}^{2d}$. The conclusion follows.

## A.4 Proof of the Lemma 2.16

The two following results are consequences of Proposition A.1:

**Corollary A.21** *Let $g \in \mathscr{C}_{pol}^\infty(\mathbb{R}^{2d})$ such that $\int_{\mathbb{R}^{2d}} g(q,p)\rho(q,p)dqdp = 0$, then there exists a unique function $\mu \in \mathscr{C}_{pol}^\infty(\mathbb{R}^{2d})$ such that*

$$L\mu = g \quad and \quad \int_{\mathbb{R}^{2d}} \mu(q,p)\rho(q,p)dqdp = 0.$$

*Proof* It is know that the unique solution $\mu$ of $L\mu = g$ which verifies $\int_{\mathbb{R}^{2d}} \mu(q,p)\rho(q,p)dqdp = 0$ is defined by $\mu(q,p) = \int_0^{+\infty} \mathbb{E}(g(q_q(t), p_p(t)))dt$. Then the regularity of $\mu$ is a consequence of Proposition A.1.

**Remark A.22** $L$ and its formal adjoint in $L^2(\rho)$ have the same behavior. Indeed, we have for any function $\phi \in C^\infty(\mathbb{R}^{2d})$ and $(q,p) \in \mathbb{R}^{2d}$

$$L^*(q,p;\partial_p,\partial_q)\phi(q,p) = -\langle p, \partial_q \phi(q,p) \rangle + \langle \partial_q V, \partial_p \phi(q,p) \rangle - \gamma \langle p, \partial_p \phi(q,p) \rangle$$
$$+ \frac{\sigma^2}{2} \sum_{i=1}^d \frac{\partial^2}{\partial_{p_i} \partial_{p_i}} \phi(q,p)$$
$$= L(q,-p;\partial_p,\partial_q)\phi(q,-p).$$

Then we have the Lemma 2.16:
Let $g \in \mathscr{C}_{pol}^\infty(\mathbb{R}^{2d})$ such that $\int_{\mathbb{R}^{2d}} g(q,p)\rho(q,p)dqdp = 0$. Then there exists a unique function $\mu \in \mathscr{C}_{pol}^\infty(\mathbb{R}^{2d})$ such that

$$L^*\mu = g \quad and \quad \int_{\mathbb{R}^{2d}} \mu(q,p)\rho(q,p)dqdp = 0,$$

where $L^*$ is the adjoint of $L$ with respect to the $L^2(\rho)$ product.

**Acknowledgements** The author is glad to thank Arnaud Debussche and Erwan Faou for their comments during the preparation of this article.